\documentclass[12pt]{amsart}
\usepackage[utf8]{inputenc}
\usepackage[margin=1in]{geometry}
\pdfoutput=1

\usepackage{bbm}
\usepackage{tikz}
\usetikzlibrary{arrows}
\usepackage{tikz-cd}
\usepackage{pgffor}
\usepackage{amssymb,amsmath,amsfonts,mathrsfs,amsthm}
\usepackage{enumitem}
\usepackage{adjustbox}
\usepackage{hyperref}
\usepackage{relsize}
\usepackage{graphicx}
\usepackage[font={small}]{caption}
\usepackage{float}
\usepackage{subcaption}
\usepackage{changepage}
\usepackage[toc]{appendix}
\usepackage{pgfplots}
\pgfplotsset{compat=1.7}
\usepackage{adjustbox}
\usepackage{varwidth}
\usepackage{chngcntr}
\usepackage{standalone}
\usepackage{tcolorbox}
\usepackage{mathtools}
\usepackage{rotating}
\usepackage{wasysym}

\usetikzlibrary{calc,knots}
\usetikzlibrary{patterns}
\usetikzlibrary{patterns.meta}
\usetikzlibrary{arrows.meta,calc,decorations.markings,math,arrows.meta, decorations.markings}
\usetikzlibrary{arrows.meta,calc,decorations.markings,math,arrows.meta}

\newcommand{\RR}{\mathbb{R}}

\newcommand{\Z}{\mathbb{Z}}
\newcommand{\C}{\mathbb{C}}

\newcommand{\A}{\mathbb{A}}

\newcommand{\id}{\operatorname{id}}

\renewcommand{\l}{\ell}

\renewcommand{\sl}{\mathfrak{sl}}
\newcommand{\gl}{\mathfrak{gl}}

\renewcommand{\d}{\partial}

\newcommand{\til}[1]{\widetilde{#1}}
\renewcommand{\b}[1]{\overline{#1}}

\newcommand{\gNmod}{{\rm -}\operatorname{gmod}}

\newcommand{\qdeg}{\operatorname{qdeg}}
\newcommand{\adeg}{\operatorname{adeg}}

\newcommand{\eps}{\varepsilon}

\renewcommand{\P}{\mathcal{P}} 
\newcommand{\Fr}{\operatorname{Fr}} 

\newcommand{\anch}{\mathrm{an}}

\newcommand{\brak}[1]{\ensuremath{\left\langle #1\right\rangle}}
\newcommand{\adm}{\operatorname{adm}}

\newcommand{\h}[1]{\widehat{#1}}

\newcommand{\OP}{\Pi}

\newcommand{\NN}{[N]}

\newcommand{\less}[2]{[{#1}\leq{#2}]}

\newcommand{\lessstrict}[2]{[{#1}<{#2}]}
\newcommand{\aFoam}{ \mathbf{AFoam}}

\renewcommand{\L}{\mathcal{L}}
\newcommand{\thick}{\operatorname{th}}
\newcommand{\Web}{\operatorname{Web}}
\newcommand{\RN}{R_N}
\newcommand{\RNsym}{R_N^{\operatorname{sym}}}

\newtheorem{theorem}{Theorem}[section]
\newtheorem{lemma}[theorem]{Lemma}
\newtheorem{proposition}[theorem]{Proposition}

    \renewcommand{\sl}{\mathfrak{sl}}

\theoremstyle{definition}

\theoremstyle{definition}

\theoremstyle{remark}
\newtheorem{example}[theorem]{Example}

\theoremstyle{remark}
\newtheorem{remark}[theorem]{Remark}

\theoremstyle{definition}
\newtheorem{definition}[theorem]{Definition}

\theoremstyle{remark}
\newtheorem{rmk/}{Remark}
\theoremstyle{definition}
\newtheorem{case/}{Case}

\usetikzlibrary{calc}

\title{Anchored $\mathfrak{gl}_N$ foams and  annular Khovanov-Rozansky homology}

\begin{document}

\begin{abstract}
    We introduce equivariant  $\mathfrak{gl}_N$ homology for links in the thickened annulus via foam evaluation. 
\end{abstract}

\author[R. Akhmechet]{Rostislav Akhmechet}
\address{Department of Mathematics, Columbia University, New York, NY 10027}
\email{\href{mailto:akhmechet@math.columbia.edu}{akhmechet@math.columbia.edu}}
\date{}

\maketitle
\tableofcontents 
\section{Introduction}

The Reshetikhin-Turaev invariant \cite{RT-link-invariant} of a link in $3$-space whose components are colored by quantum exterior powers  of the fundamental quantum $\sl_N$ representation can be computed by resolving crossings in a link diagram, resulting in a $\Z[q^{\pm 1}]$-linear combination of certain trivalent graphs called Murakami-Ohtsuki-Yamada (MOY) \emph{webs}, and then evaluating each web $\Gamma$ to a Laurent polynomial $\brak{\Gamma}_N \in \Z_{\geq 0} [q^{\pm 1}]$ \cite{MOY}. There are now many approaches to categorifying this link polynomial, beginning with Khovanov homology \cite{Kh} in the case of $N=2$ and the Jones polynomial, Khovanov's $\sl_3$ homology \cite{Khsl3}, and  Khovanov-Rozansky homology \cite{KRoz} for general $N$ and when each color is the fundamental representation. 

We focus on the approach established by Robert and Wagner \cite{RWfoamev}. To a MOY web $\Gamma$ they assign a module  $\brak{\Gamma}_{\rm RW}$, called the \emph{state space},  which categorifies the MOY calculus in the sense that $\brak{\Gamma}_{\rm RW}$ is a free graded module of graded rank equal to $\brak{\Gamma}_N$ \cite[Theorem 3.30]{RWfoamev}. Given a link diagram, one can form a chain complex using a cube of resolutions, where each vertex is decorated by the state space of a resolution web of the link diagram. This yields a strictly functorial homology theory of links in $\RR^3$, proven by Ehrig-Tubbenhauer-Wedrich \cite{ETW}. Other approaches to categorifying the exterior colored $\sl_N$ polynomial include \cite{MSVfoams, Yonezawa,Wu, Wu_equivariant, QRfoams}. Webster \cite{Webster} categorified the Reshetikhin-Turaev invariant for any simple Lie algebra $\mathfrak{g}$ and components colored by any representation.

The key construction is that of the state spaces $\brak{\Gamma}_{\rm RW}$. In \cite{RWfoamev} these are built combinatorially: one first considers \emph{foams}, which are certain $2$-dimensional CW complexes, viewed as cobordisms between webs. Robert and Wagner introduced a combinatorial formula to evaluate closed  foams. The state spaces are then defined by applying the  universal construction, introduced in \cite{BHMV}, to this closed foam evaluation. It is more natural to view the theory as $\gl_N$ rather than $\sl_N$ link homology, so in this paper we will refer to $\gl_N$ webs and foams. Foam evaluation  has seen many interesting developments in recent years \cite{KRfoamev, RWsymm, RWAlexander, Khovanov-Kitchloo,  BPRW, QRSWsl2}.

Link homology built using Robert-Wagner foam evaluation is naturally \emph{equivariant}, in the sense that an $a$-colored unknot is assigned the $U(N)$-equivariant cohomology of the Grassmannian  of $a$-dimensional subspaces in $\C^N$. The ground ring is the $U(N)$-equivariant cohomology of a point, which can be identified with $\RNsym = \Z[x_1, \ldots, x_N]^{S_N}$, the graded ring of symmetric polynomials with each $x_i$ in degree $2$. State spaces of webs and homology groups of links are graded $\RNsym$-modules. One may specialize each $x_i$ to $0$ to obtain graded abelian groups. Deformations were studied in \cite{RWdeformations}, and generic deformations are important for establishing strict functoriality \cite{ETW}. 

Let us briefly discuss the case $N=2$. We have $R_2^{\operatorname{sym}} = \Z[E_1, E_2]$, where $E_1 = x_1 + x_2$ and $E_2 = x_1 x_2$ are the elementary symmetric polynomials. The  state space of a $1$-colored circle is a Frobenius algebra given by
\[
\brak{\mathlarger{\fullmoon}_1}_{\rm RW} \cong \Z[E_1, E_2, X]/(X^2 - E_1 X + E_2),
\]
which, upon renaming variables $E_1 =  h, E_2 = -t$, is precisely the Frobenius algebra denoted $\mathcal{F}_5$ in \cite{KhFrobext}. Specializing the variables $x_i$ to $0$ recovers the $\gl_2$ link homology originally defined by Blanchet \cite{Blanchetgl2}. Generic specializations of $x_1, x_2$ lead to an analogue of the Lee deformation \cite{Lee}.

This paper focuses on links in the thickened annulus. Beginning with work of Asaeda-Przytycki-Sikora \cite{APS},  Khovanov homology  for links in the thickened annulus has been a subject of  interest \cite{Roberts,GLWsl2, GN, GLWsinvt, BG, BPW}. The so-called \emph{annular} Khovanov (or annular APS) homology is triply graded: in addition to the homological and quantum gradings, there is a third  \emph{annular} grading. Annular versions of Khovanov-Rozansky  homology  were developed in \cite{QR, QRS}.

An equivariant version of annular $\sl_2$ homology was defined in previous work of the author \cite{Akh}; this did not use foam evaluation but rather a filtration in the spirit of \cite{Roberts, GLWsl2}. Then in joint work with Khovanov \cite{AkhKh} a foam evaluation approach to equivariant annular $\sl_2$ and $\sl_3$ link homology was introduced. The annular $\sl_2$ homology defined via foam evaluation was shown to be isomorphic to the one constructed in \cite{Akh}. A key phenomenon appearing in the annular setting is the need to extend the ground ring from symmetric polynomials to all polynomials.

The main idea is to identify the interior of an annulus with a punctured plane  and the interior of a thickened annulus with the complement of a line $\L$ in $\RR^3$. We choose for convenience the punctured plane $\P = \RR^2\setminus \{(0,0)\}$ and the line $\L$ to be the $z$-axis. For a web $\Gamma$ in $\P$, its state space is spanned by foams in $\RR^2 \times (-\infty, 0]$ whose boundary is $\Gamma$. Foams may intersect the line $\L$ generically, and these intersection points, called \emph{anchor points},  carry  decorations which contribute additional factors to the foam evaluation. 

In the present paper we extend the construction in \cite{AkhKh} to anchored $\gl_N$ foams. An anchored $\gl_N$ foam may intersect $\L$ transversely in the interior of its facets. Each anchor point $p$ carries a label consisting of a subset of $\{1, \ldots, N\}$ of cardinality equal to the thickness of the facet on which $p$ lies. Anchored $\gl_N$ foam evaluation takes values in the ring $\Z[x_1, \ldots, x_N]$ rather than in its subring of symmetric polynomials; consequently, state spaces of annular webs and  homology of links in the thickened annulus are modules over this larger ring. 

In addition to the  quantum grading, state spaces of annular $\gl_N$ webs also carry a $\Z^N$ \emph{annular} grading  coming from labels of anchor points, analogous to the annular gradings in \cite{AkhKh}. The annular $\gl_N$ link homology defined in this paper is thus $\Z \oplus \Z^N$-graded.  The $\Z^N$-grading is expected from  work of Queffelec-Rose \cite{QR}, which introduces (non-equivariant) annular $\sl_N$ homology that admits an action of $\sl_N$. The annular grading then corresponds to the weight space decomposition. The work of Queffelec-Rose generalizes results of Grigsby-Licata-Wehrli \cite{GLWsl2}, who showed that annular Khovanov ($\sl_2$) homology carries an action of $\sl_2$ in which the annular grading corresponds to the weight space grading. 

As noted in \cite{AkhKh}, it is not clear how to extend actions of $\sl_2$ and $\sl_3$ to the equivariant annular theories, and the same is true in the present paper in the context of $\gl_N$ foams. In a related but different direction,  Qi-Robert-Sussan-Wagner \cite{QRSWsl2} established an action of a subalgebra of the Witt algebra on $\gl_N$ foams and on state spaces of $\gl_N$ webs. As a consequence, they obtain an action of $\sl_2$ on $\gl_N$ state spaces.
In future work we will extend this action to the annular setting.

This paper is organized as follows. In Section \ref{sec:Webs and foams} we review $\gl_N$ webs, foams, and Robert-Wagner closed foam evaluation. Anchored $\gl_N$ foams and their evaluation are defined in Section \ref{sec:anchored gln foam evaluation};  Proposition \ref{prop:no denominators} establishes that the evaluation  is a polynomial. In Section \ref{sec:anchored gln local relations} we establish skein relations on anchored foams. Section \ref{sec:state spaces annular gln} discusses how to build state spaces of annular webs using anchored foam evaluation and universal construction. The main result is Theorem \ref{thm:anchored gln state space}, which identifies state spaces of annular webs. In Section \ref{sec:equivariant annular gln link homology} we discuss how to obtain equivariant colored annular homology as well as functoriality properties. The content of this paper is largely taken from the author's Ph.D. thesis \cite{Akhthesis}.

\vskip3ex

\textbf{Acknowledgments:} I would like to thank Mikhail Khovanov for encouraging the pursuit of this project and for many enlightening discussions on foams. I was partially supported by NSF Grant DMS-2105467 and NSF RTG grant DMS-1839968 while working on this project. 

\section{Webs and foams}
\label{sec:Webs and foams}

\subsection{MOY webs}

Fix an integer $N\geq 1$. In this section we discuss Murakami-Ohtsuki-Yamada (MOY) webs \cite{MOY}. These are certain planar trivalent graphs which can be used to compute the exterior colored Reshetikhin-Turaev $\sl_N$ (or $\gl_N$) link polynomial and which give planar diagrammatics describing morphism spaces between tensor products of exterior powers of the defining $U_q(\sl_N)$ representation.

\begin{definition}
For $n\in \Z$, the quantum integer $[n]_q$ is defined to be 
\[
[n]_q= \frac{q^n - q^{-n}}{q-q^{-1}} = q^{n-1} + q^{n-3}  + \cdots + q^{3-n} + q^{1-n} \in \Z_{\geq 0}[q,q^{-1}].
\]
If $k>0$, set 
\[
\genfrac[]{0pt}{0}{n}{k}_q = \frac{[n]_q[n-1]_q\cdots [n-k+1]_q}{[k]_q[k-1]_q\cdots[1]_q},
\]
and define $\genfrac[]{0pt}{0}{n}{0}_q=1$. 
\end{definition}

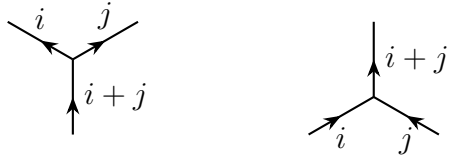
\begin{figure}
\centering
{\begin{tikzpicture}[thick,decoration={
    markings,
    mark=at position 0.5 with {\arrow{Stealth}}}
    ] 

\draw[postaction={decorate}] (0,-1) -- (0,0) node [midway, right] {$i+j$};
\draw[postaction={decorate}] (0,0) -- ({cos(30)},{sin(30)}) node[midway,above] {$j$};
\draw[postaction={decorate}] (0,0) -- (-{cos(30)}, {sin(30)}) node[midway,above] {$i$};

\begin{scope}[shift={(4,{-sin(30)})}]
\draw[postaction={decorate}] ({-cos(30)},{-sin(30)}) -- (0,0) node [midway,below] {$i$};
\draw[postaction={decorate}] ({cos(30)},{-sin(30)})-- (0,0) node[midway,below] {$j$};
\draw[postaction={decorate}] (0,0) -- (0,1) node[midway,right] {$i+j$};
\end{scope}

\end{tikzpicture}
}
\caption{The flow condition near each trivalent vertex in a $\gl_N$ web.}\label{fig:web flow condition}
\end{figure}

\begin{definition}
A \emph{$\gl_N$ web} (also called a \emph{MOY web}) is a embedded trivalent graph $\Gamma \subset \RR^2$, which may also contain closed loops with no vertices. Moreover, edges and loops of $\Gamma$ are oriented and carry weights in $\{0,\ldots, N\}$, called the \emph{thickness} of the edge, such that the flow condition shown in Figure \ref{fig:web flow condition} is satisfied at each vertex. Let $\Web_N$ denote the set of planar isotopy classes of $\gl_N$ webs. We will often simply write \emph{web} instead of $\gl_N$ web. 
\end{definition}

\begin{figure}
\centering
\subcaptionbox{
\label{MOYa}
 }[.4\linewidth]
{ \begin{tikzpicture}[decoration={
    markings,
    mark=at position 0.5 with {\arrow{Stealth}}}
    ] 

\draw[postaction={decorate}] (0,0) arc (-90:90:.75) node [midway,right] {$i$};
\draw (0,1.5) arc (90:270:.75);

\node at (2,.75) {$=\ \ \genfrac[]{0pt}{0}{N}{i}_q$};
\end{tikzpicture}
}
\subcaptionbox{
\label{MOYb}
 }[.4\linewidth]
{\begin{tikzpicture}[yscale=1.5,decoration={
    markings,
    mark=at position 0.7 with {\arrow{Stealth}}}
    ] 
\draw[postaction={decorate}] (0,-.75) -- (0,0) node [pos=.2, right] {$\mathsmaller{i+j+k}$};

\draw[postaction={decorate}] (0,0) -- ({cos(30)/2},{sin(30)/2}) node[midway,below] {$\ \ \ \mathsmaller{j+k}$};

\draw[postaction={decorate}] ({cos(30)/2},{sin(30)/2})--++ (-{cos(30)/2},{sin(30)/2}) node[midway,above] {$ \mathsmaller{j}$};

\draw[postaction={decorate}] ({cos(30)/2},{sin(30)/2})--++ ({cos(30)/2},{sin(30)/2}) node[midway,above] {$ \mathsmaller{k}$};

\draw[postaction={decorate}] (0,0) -- (-{cos(30)}, {sin(30)}) node[midway,above] {$\mathsmaller{i}$};


\begin{scope}[shift={(3,0)}]
\draw[postaction={decorate}] (0,-.75) -- (0,0) node [pos=.2, right] {$\mathsmaller{i+j+k}$};

\draw[postaction={decorate}] (0,0) -- ({cos(30)},{sin(30)}) node[midway,above] {$\mathsmaller{k}$};

\draw[postaction={decorate}] (0,0) -- (-{cos(30)/2}, {sin(30)/2}) node[midway,below] {$\mathsmaller{i+j}\ \ \ $};

\draw[postaction={decorate}] (-{cos(30)/2}, {sin(30)/2})--++ ({cos(30)/2},{sin(30)/2}) node[midway,above] {$ \mathsmaller{j}$};

\draw[postaction={decorate}] (-{cos(30)/2}, {sin(30)/2})--++ (-{cos(30)/2},{sin(30)/2}) node[midway,above] {$ \mathsmaller{i}$};

\end{scope}

\node at (1.5,0) {$=$};
\end{tikzpicture}
}\\ \vskip 1ex
\subcaptionbox{
\label{MOYc}
 }[.4\linewidth]
{\begin{tikzpicture}[yscale = .8, decoration={
    markings,
    mark=at position 0.5 with {\arrow{Stealth}}}
    ] 

\draw[postaction={decorate}] (0,0) -- (0,1) node[midway,right] {$i+j$};

\draw[postaction={decorate}] (0,-2) -- (0,-1) node[midway,right] {$i+j$};
    
\draw[postaction ={decorate}] (0,-1) to [bend right=80] node[midway, right] {$j$} (0,0)  ;

\draw[postaction ={decorate}] (0,-1) to [bend left=80] node[midway, left] {$i$} (0,0)  ;


\begin{scope}[shift={(3.5,0)}]
\node at (-1,-.5) {$\genfrac[]{0pt}{0}{i+j}{i}_q$}; 

\draw[postaction={decorate}] (0,-2) -- (0,1) node[midway,right] {$i+j$};
\end{scope}

\node at (1.5,-.5) {$=$};
\end{tikzpicture}
}
\subcaptionbox{
\label{MOYd}
 }[.4\linewidth]
{\begin{tikzpicture}[yscale = .8,decoration={
    markings,
    mark=at position 0.5 with {\arrow{Stealth}}}
    ] 

\draw[postaction={decorate}] (0,0) -- (0,1) node[midway,right] {$i$};

\draw[postaction={decorate}] (0,-2) -- (0,-1) node[midway,right] {$i$};
    
\draw[postaction ={decorate}] (0,0) to [bend left=80] node[midway, right] {$j$} (0,-1)  ;

\draw[postaction ={decorate}] (0,-1) to [bend left=80] node[midway, left] {$i+j$} (0,0)  ;


\begin{scope}[shift={(3.7,0)}]
\node at (-1,-.5) {$\genfrac[]{0pt}{0}{N-i}{j}_q $}; 

\draw[postaction={decorate}] (0,-2) -- (0,1) node[midway,right] {$i$};
\end{scope}

\node at (1.5,-.5) {$=$};
\end{tikzpicture}
}\\ \vskip 1ex
\subcaptionbox{
\label{MOYe}
 }[.8\linewidth]
{\begin{tikzpicture}[yscale = .8,decoration={
    markings,
    mark=at position 0.5 with {\arrow{Stealth}}}
    ]

\draw[postaction={decorate}] (0,0) --++ (0,1) node[midway, left] {$\mathsmaller{1}$}; 

\draw[postaction={decorate}] (1.5,0) --++ (-1.5,0) node[midway, above] {$\mathsmaller{i+1}$}; 

\draw[postaction={decorate}] (1.5,1) --++ (0,-1) node[midway, right] {$\mathsmaller{i}$}; 

\draw[postaction={decorate}] (0,0) --++ (0,-1) node[midway,left] {$\mathsmaller{i}$}; 

\draw[postaction={decorate}] (0,-1) --++ (1.5,0) node[midway,below] {$\mathsmaller{i+1}$};

\draw[postaction={decorate}] (1.5,-1) --++ (0,1) node[midway,right] {$\mathsmaller{1}$};

\draw[postaction={decorate}] (0,-2) --++ (0,1) node[midway,left] {$\mathsmaller{1}$}; 

\draw[postaction={decorate}] (1.5,-1) --++ (0,-1) node[midway,right] {$\mathsmaller{i}$};


\begin{scope}[shift={(3.5,0)}]
\draw[postaction={decorate}] (0,-2) --++ (0,3) node[midway,left]  {$1$} ;

\draw[postaction={decorate}] (1.5,1) --++ (0,-3) node[midway,right]  {$i$} ;
\end{scope}


\begin{scope}[shift={(9,0)}]
\node at (-1.25,-.5) {$[N-i-1]_q$}; 

\draw[postaction={decorate}] (0,-2) --++ (.75,1) node[midway,left]  {$1\ $} ;

\draw[postaction={decorate}] (.75,-1) --++ (.75,-1) node[midway,right] {$\ i$};

\draw[postaction={decorate}] (.75,0) --++ (0,-1) node[midway,right] {$i-1$};

\draw[postaction={decorate}] (.75,0) --++ (-.75,1) node[midway,left] {$1\ $};

\draw[postaction={decorate}] (1.5,1) --++ (-.75,-1) node[midway,right] {$\ i$};
\end{scope}

\node at (2.5,-.5) {$=$};

\node at (6.25,-.5) {$+$};

\end{tikzpicture} 
}\\ \vskip 1ex
\subcaptionbox{
\label{MOYf}
 }[.8\linewidth]
{\begin{tikzpicture}[yscale = .8,decoration={
    markings,
    mark=at position 0.5 with {\arrow{Stealth}[scale=1]}}
    ]

\draw[postaction={decorate}] (0,0) --++ (0,1) node[midway, left] {$\mathsmaller{i}$}; 

\draw[postaction={decorate}] (0,0) --++ (1.5,0) node[midway, above] {$\mathsmaller{j}$}; 

\draw[postaction={decorate}] (1.5,0) --++ (0,1) node[midway, right] {$\mathsmaller{k}$}; 

\draw[postaction={decorate}] (0,-1) --++ (0,1) node[midway,left] {$\mathsmaller{i+j}$}; 

\draw[postaction={decorate}] (1.5,-1) --++ (-1.5,0) node[midway,below] {$\mathsmaller{i+j-1}$};

\draw[postaction={decorate}] (1.5,-1) --++ (0,1) node[midway,right] {$\mathsmaller{k-j}$};

\draw[postaction={decorate}] (0,-2) --++ (0,1) node[midway,left] {$\mathsmaller{1}$}; 

\draw[postaction={decorate}] (1.5,-2) --++ (0,1) node[midway,right] {$\mathsmaller{k+i-1}$};


\begin{scope}[shift={(5.5,0)}]
\node at (-1.25,-.5) {$\genfrac[]{0pt}{0}{k-1}{j}_q$}; 

\draw[postaction={decorate}] (0,-2) --++ (0,1.5) node[midway,left]  {$1$} ;

\draw[postaction={decorate}] (0,-.5) --++ (0,1.5) node[midway,left]  {$i$} ;

\draw[postaction={decorate}] (1.5,-2) --++ (0,1.5) node[midway,right]  {$k+i-1$} ;

\draw[postaction={decorate}] (1.5,-.5) --++ (-1.5,0) node[midway,above] {$i-1$};

\draw[postaction={decorate}] (1.5,-.5) --++ (0,1.5) node[midway,right] {$k$};
\end{scope}


\begin{scope}[shift={(11.5,0)}]
\node at (-1.25,-.5) {$\genfrac[]{0pt}{0}{k-1}{j-1}_q$}; 

\draw[postaction={decorate}] (0,-2) --++ (.75,1) node[midway,left]  {$1\ $} ;

\draw[postaction={decorate}] (1.5,-2) --++ (-.75,1) node[midway,right] {$\ k+i-1$};

\draw[postaction={decorate}] (.75,-1) --++ (0,1) node[midway,right] {$k+i$}; 

\draw[postaction={decorate}] (.75,0) --++ (-.75,1) node[midway,left] {$i\ $};

\draw[postaction={decorate}] (.75,0) --++ (.75,1) node[midway,right] {$\ k$};
\end{scope}

\node at (3,-.5) {$=$};

\node at (9,-.5) {$+$};

\end{tikzpicture} 
}\\ \vskip 1ex
\subcaptionbox{
\label{MOYg}
 }[.8\linewidth]
{\begin{tikzpicture}[yscale = .8,decoration={
    markings,
    mark=at position 0.5 with {\arrow{Stealth}}}
    ]

\draw[postaction={decorate}] (0,0) --++ (0,1) node[midway, left] {$\mathsmaller{i}$}; 

\draw[postaction={decorate}] (0,0) --++ (1.5,0) node[midway, above] {$\mathsmaller{j+k-i}$}; 

\draw[postaction={decorate}] (1.5,0) --++ (0,1) node[midway, right] {$\mathsmaller{j+l}$}; 

\draw[postaction={decorate}] (0,-1) --++ (0,1) node[midway,left] {$\mathsmaller{j+k}$}; 

\draw[postaction={decorate}] (1.5,-1) --++ (-1.5,0) node[midway,below] {$\mathsmaller{k}$};

\draw[postaction={decorate}] (1.5,-1) --++ (0,1) node[midway,right] {$\mathsmaller{i+l-k}$};

\draw[postaction={decorate}] (0,-2) --++ (0,1) node[midway,left] {$\mathsmaller{j}$}; 

\draw[postaction={decorate}] (1.5,-2) --++ (0,1) node[midway,right] {$\mathsmaller{i+l}$};


\begin{scope}[shift={(7.5,0)}]
\node at (-2.5,-.5) {$\displaystyle \sum\limits_{n=\max(0,i-j)}^j \ \genfrac[]{0pt}{0}{l}{k-n}_q $};

\draw[postaction={decorate}] (0,0) --++ (0,1) node[midway, left] {$\mathsmaller{i}$}; 

\draw[postaction={decorate}] (1.5,0) --++ (-1.5,0) node[midway, above] {$\mathsmaller{n}$}; 

\draw[postaction={decorate}] (1.5,0) --++ (0,1) node[midway, right] {$\mathsmaller{j+l}$}; 

\draw[postaction={decorate}] (0,-1) --++ (0,1) node[midway,left] {$\mathsmaller{i-n}$}; 

\draw[postaction={decorate}] (0,-1) --++ (1.5,0) node[midway,below] {$\mathsmaller{j-i+n}$};

\draw[postaction={decorate}] (1.5,-1) --++ (0,1) node[midway,right] {$\mathsmaller{j+l+n}$};

\draw[postaction={decorate}] (0,-2) --++ (0,1) node[midway,left] {$\mathsmaller{j}$}; 

\draw[postaction={decorate}] (1.5,-2) --++ (0,1) node[midway,right] {$\mathsmaller{i+l}$};
\end{scope}


\node at (3,-.5) {$=$};

\end{tikzpicture} 
}
\caption{The MOY relations.}\label{fig:MOY relations}
\end{figure}

\begin{theorem}[{\cite{MOY}, see also \cite[Theorem 2.4]{Wu}}]
There is a unique function $\brak{-}_N : \Web_N \to \Z_{\geq 0}[q^{\pm 1}]$ which factors through the local relations in Figure \ref{fig:MOY relations}.
\end{theorem}

Indeed, the evaluation $\brak{\Gamma}_N$ for any web $\Gamma$ can be  computed using the relations in Figure \ref{fig:MOY relations}. Alternatively, this evaluation may be computed by applying a state-sum formula as in \cite{MOY}.

\subsection{Robert-Wagner foam evaluation}
\label{sec:RW foam evaluation}

This section gives an overview of Robert-Wagner foam evaluation, which is used to build state spaces of webs that categorify the MOY calculus. Later in Section \ref{sec:anchored gln} we extend foam evaluation to the annular setting.

\begin{figure}
\centering
\subcaptionbox{A binding, where three facets meet. \label{fig:foam_orientation}}[.4\linewidth]
{\begin{tikzpicture}[decoration={
    markings,
    mark=at position 0.6 with {\arrow{Stealth[scale=1]}}}
    ] 
\draw (0,0) -- (-1,1) -- (1,3) -- (2,2);

\draw (0,0) -- (2.5,1) -- (4.5,3) -- (2,2);

\draw (0,0) -- (0,-1.5) -- (2,.5) -- (2,.75);

\draw[densely dashed] (2,2) -- (2,.75);

\draw[very thick, postaction=decorate] (0,0) -- (2,2);

\draw[->] (.5,1.75) arc (110:430:.2 and .2);

\draw[<-] (2.5,1.75) arc (110:430:.2 and .2);

\draw[->] (1,.15) arc (110:430:.2 and .2);

\node at (1,2.1) {$\mathsmaller{a}$}; 
\node at (3,2.1) {$\mathsmaller{b}$}; 

\node at (.5,-.5) {$\mathsmaller{a+b}$};

\end{tikzpicture}
}
\subcaptionbox{A singular vertex, where six facets meet. \label{fig:singular vertex}}[.5\linewidth]
{\begin{tikzpicture}[scale=.9,decoration={
    markings,
    mark=at position 0.6 with {\arrow{Stealth}}}
    ] 
\draw[fill=black] (0,0) circle [radius=0.1];

\draw[very thick, postaction={decorate}] (-3,1) -- (-2,{2/3});
\draw[very thick, densely dashed] (-2,{2/3}) -- (0,0);
\draw[very thick,postaction={decorate}](0,0) -- (3,-1);
\draw[very thick, postaction={decorate}] (-2,-1) -- (0,0);

\draw[very thick, postaction={decorate}] (0,0) -- (2,1);

\draw(-2,-1) -- (1,-2);
\draw (1,-2) -- (3,-1);

\draw (3,-1) -- (5,0);
\draw(5,0) -- (2,1);

\draw(-5,0) -- (-2,-1);
\draw (-5,0) -- (-3,1);

\draw (-3,1) -- (-1,2);
\draw (-1,2) -- (-1+.6, {2- .6/3});
\draw[dashed] (-1+.6, {2- .6/3}) -- (2,1);

\draw(-2,-1) -- (-2,1);
\draw(2,1) -- (2,3);
\draw (-2,1) -- (2,3);

\draw(3,-1) -- (3,-3);
\draw(3,-3) -- (-3,-1);
\draw(-3,-1) -- (-3,-.65);
\draw[dashed] (-3,-.65) -- (-3,1);

\node at (-1,1.75) {$\mathsmaller{a}$}; 
\node at (-3.75,0) {$\mathsmaller{b+c}$};

\node at (1.75,2.5) {$\mathsmaller{b}$}; 

\node at (3.75,0) {$\mathsmaller{a+b}$}; 

\node at (1,-1.5) {$\mathsmaller{c}$};

\node at ((2.2,-2.1) {$\mathsmaller{a+b+c}$};
\end{tikzpicture}
}
\caption{A summary of the thickness and orientation conventions at singular points.}\label{fig:singular points}
\end{figure}
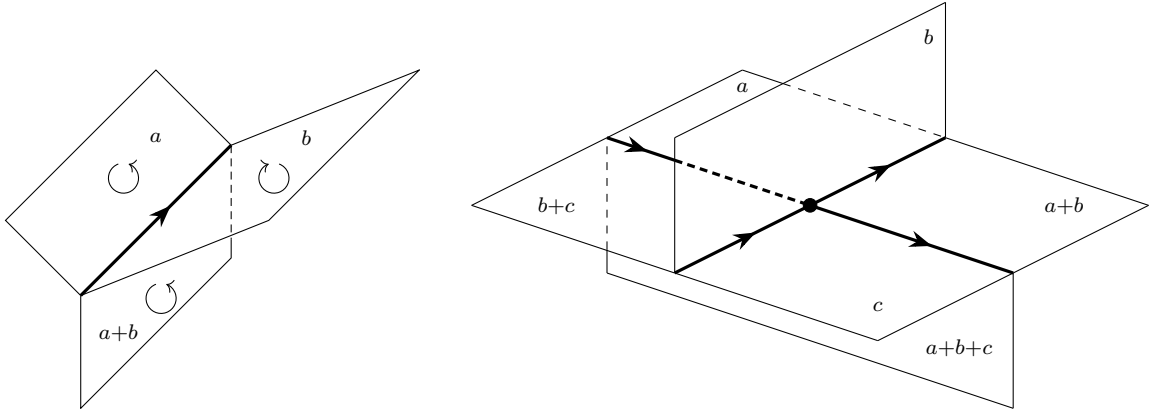

\begin{definition}
A (closed) \emph{$\gl_N$ foam} $F$ is a compact, PL, 2-dimensional CW complex with the following properties and data.

\begin{itemize}
    \item  Every point of $F$ has neighborhood homeomorphic to either a disk, the letter $Y$ times an interval, or the cone on the 1-skeleton of a tetrahedron; the latter two types of points, illustrated in Figure \ref{fig:singular points}, are called \emph{singular} points. The \emph{singular graph of $F$}, denoted $s(F)$, consists of all singular points. Note that $s(F)$ may contain closed loops with no vertices. Edges and loops in $s(F)$ are called \emph{bindings} or \emph{seams}; moreover, the edges are called \emph{interval bindings}. A \emph{singular vertex} of $F$ is a vertex of $s(F)$, shown in Figure \ref{fig:singular vertex}. 

\item A \emph{facet} of $F$ is a connected component of $F\setminus s(F)$; the set of facets is denoted $f(F)$. Each facet $f\in f(F)$ of $F$ carries a label in $\{0,\ldots, N\}$ called its \emph{thickness} and denoted $\thick(f)$. 

\item The singular graph and facets of $F$ must be oriented, and moreover these orientations must be compatible with thicknesses in the following way. Three facets meeting along a binding must have thickness $a,b,$ and $a+b$. The orientation of the binding must agree with the orientations of the thickness $a$ and $b$ facets and disagree with the orientation of thickness $a+b$ facet. See Figure \ref{fig:foam_orientation} for a summary.

\item Foams moreover carry \emph{decorations}, consisting of a homogeneous symmetric polynomial $P_f$ in $\thick(f)$ variables for each facet $f$. We will indicate these polynomial decorations by an arrow pointing from a polynomial to the facet it decorates, see Figure \ref{fig:theta foam unanchored}. When depicting foams, a facet with no such decoration appearing indicates that it is decorated by the constant polynomial $1$. 

\item Finally, foams considered in the present paper will always be PL embedded in $\RR^3$. As for webs, we will often write \emph{foam} instead of $\gl_N$ foam. 
\end{itemize}
\end{definition}

\begin{remark}
One may consider abstract foams (not embedded in $\RR^3$). In this case a cyclic order of the three facets at each interval binding must be fixed. In our situation, when all foams are embedded, the cyclic order is determined by the left-hand rule. 
\end{remark}

\begin{example}
Consider the so-called \emph{theta foam}  $F$ shown in Figure \ref{fig:theta foam unanchored}. It consists of three disks with thicknesses $a,b$, and $a+b$ glued along their common boundary, with symmetric polynomial labels $P_1, P_2$, and $P_3$, respectively. The orientations on facets can be determined from the orientation of the singular circle.

\begin{figure}
    \centering
    \begin{tikzpicture}[scale=1.7,decoration={
    markings,
    mark=at position 0.75 with {\arrow{Stealth[scale=2]}}}]

\draw (0,0) circle [radius=1];

\draw[densely dashed,fill=gray!40] (0,0) ellipse (1 and .3);
\draw[densely dashed] (1,0) arc(0:180:1 and .3);
\draw[postaction={decorate}] (1,0) arc(0:-180:1 and .3);

\node(dt) at (.4,.6) {$b$};
\node(dm) at (.4,0) {$a+b$};
\node(db) at (.4,-.6) {$a$};

\node(lt) at (-1.5,.6) {$P_2$};
\node(lm) at (-1.5,0) {$P_3$};
\node(lb) at (-1.5,-.6) {$P_1$};

\draw[thick, ->, blue] (lt) -- (-.6,.6);
\draw[thick, ->, blue] (lm) -- (-.6,0);
\draw[thick, ->, blue] (lb) -- (-.6,-.6);
\end{tikzpicture}
    \caption{A decorated theta foam.}
    \label{fig:theta foam unanchored} 
\end{figure}
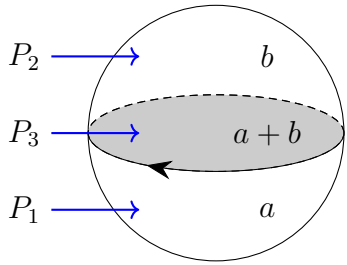

\end{example}

\begin{definition}
\label{def:qdeg of gln foam}
Let $F$ be a $\gl_N$ foam.  Consider the following quantities. 
\begin{itemize}
    \item For a facet $f$ of thickness $a$, set $d(f) = a(N-a)\chi(f)$.
    \item For an interval binding $i$ which meets facets of thicknesses $a, b$, and $a+b$, set $d(i) = ab + (a+b)(N-a-b)$. 
    \item For a singular vertex $v$ which meets facets of thicknesses $a,b, c, a+b, a+c, b+c$, and $a+b+c$, set $d(v) = ab+ac+bc + (N-a-b-c)(a+b+c)$.
\end{itemize}
Define the degree of a foam $F$ to be
\[
\deg(F) = -\sum_{f\in f(F)}  d(f)
+
\sum_{\substack{i {\rm\  an } \\ {\rm interval \ binding}}} d(i) 
-
\sum_{\substack{v {\rm\  a } \\ {\rm singular\  vertex}}} d(i) +
\sum_{f\in f(F)} \deg(P_f)
\]
\end{definition}

We record the following definition here, although it will not be used until Section \ref{sec:state spaces annular gln}.

\begin{definition}
\label{def:gln foam with boundary}
Let $\Gamma_0, \Gamma_1$ be $\gl_N$ webs. A \emph{$\gl_N$ foam with boundary $(\Gamma_0, \Gamma_1)$} consists of an intersection of a $\gl_N$ foam $F$ with $\RR^2\times [0,1]$ such that, for some $\eps>0$, $F\cap (\RR^2\times [1-\eps, 1]) = \Gamma_1 \times [1-\eps, \eps]$ and $F\cap (\RR^2\times [0,\eps]) = (-\Gamma_0) \times [0,\eps]$,  where $-\Gamma_0$ denotes the embedded graph $\Gamma_0$ with orientation reversed. 
\end{definition}

We introduce the following notation. 
\begin{itemize}
    \item Let $\NN= \{1, 2, \ldots, N\}$. 
    \item For a set $A$,  let $2^A$ denote its powerset.
    \item We write $\# A$ to denote the cardinality of a finite set $A$.     
\end{itemize}

\begin{definition}
\label{def:gln foam coloring}
An \emph{coloring} of a foam $F$ is a function $c: f(F) \to 2^{\NN}$ such that $\#c(f) = \thick(f)$ for all $f\in f(F)$. Moreover, if three facets $f_1, f_2, f_3$, of thicknesses $a,b$ and $a+b$ respectively, meet at a singular edge, then we must have $c(f_1) \sqcup c(f_2) = c(f_3)$. This condition is illustrated in Figure \ref{fig:coloring}. Let $\adm(F)$ denote the set of  colorings of $F$.

For $i\in \NN$, we say a facet $f\in f(F)$ is \emph{colored $i$} (according to $c$) if $i\in c(f)$; more generally, for $I\subset \NN$, we say $f$ is colored $I$ if $I\subset c(f)$. 
\end{definition}

Now fix a foam $F$ and a coloring $c$. 
\begin{itemize}
    \item For $1\leq i \leq N$, let $F_i(c)$ denote the union of all $i$-colored facets of $F$. 
    \item For $1\leq i\neq j \leq N$, let $F_{ij}(c)$ denote the union of facets colored either $i$ or $j$ but not both. 
    
\item Let $1\leq i<j \leq N$ and let $\beta$ be a binding joining facets $f_1, f_2, f_3$ with $f_1$ colored $i$, $f_2$ colored $j$, and $f_3$ colored $\{i,j\}$. We say $\beta$ is \emph{positive} with respect to $(i,j)$ if, according to the left-hand rule, the cyclic ordering is $(f_1, f_2, f_3)$; otherwise we say $\beta$ is negative. This is summarized in Figure \ref{fig:positive binding}. Let 
\[
\theta^+_{ij}(c)
\]
denote the number of positive $(i.j)$ bindings.

\item Recall that facets of foams are decorated by symmetric polynomials. If $A\subset \NN$ and $P \in \Z[y_1, \ldots, y_m]$ is a symmetric polynomial in $m=\# A$ variables, let $P(A) := P((x_a)_{a\in A}) \in \Z[x_1, \ldots, x_N]$. 

\end{itemize}

\begin{figure}
\centering
\subcaptionbox{The condition on a foam coloring $c$ around each singular edge.
\label{fig:coloring}}[.55\linewidth]
{\begin{tikzpicture}[xscale=1,yscale=1,decoration={
    markings,
    mark=at position 0.6 with {\arrow{Stealth[scale=1]}}}
    ] 
\draw (0,0) -- (-1,1) -- (1,3) -- (2,2);

\draw (0,0) -- (2.5,1) -- (4.5,3) -- (2,2);

\draw (0,0) -- (0,-1.5) -- (2,.5) -- (2,.75);

\draw[densely dashed] (2,2) -- (2,.75);

\draw[very thick, postaction=decorate] (0,0) -- (2,2);

\node at (1,2.1) {${a}$}; 
\node at (3,2.1) {${b}$}; 

\node at (.5,-.3) {${a+b}$};

\draw[->,blue] (-1,1.5) -- (.5,1.5); 

\draw[->,blue] (4,1.5) -- (2.5,1.5);

\draw[->,blue] (2.5,0) -- (1,0);

\node[left] at (-1,1.5) {$c(f_1)$};

\node[right] at (4,1.5) {$c(f_2)$};

\node[right] at (2.5,0) {$c(f_3) = c(f_1) \cup c(f_2)$};

\end{tikzpicture}
}
\subcaptionbox{A positive $(i,j)$ binding, where $1\leq i < j \leq N$.
\label{fig:positive binding}}[.4\linewidth]
{\begin{tikzpicture}[decoration={
    markings,
    mark=at position 0.6 with {\arrow{Stealth[scale=1]}}}
    ] 
\draw (0,0) -- (-1,1) -- (1,3) -- (2,2);

\draw (0,0) -- (2.5,1) -- (4.5,3) -- (2,2);

\draw (0,0) -- (0,-1.5) -- (2,.5) -- (2,.75);

\draw[densely dashed] (2,2) -- (2,.75);

\draw[very thick, postaction=decorate] (0,0) -- (2,2);

\node at (1,2) {$\{j\}$}; 
\node at (3,2) {$\{i\}$}; 

\node at (.5,-.4) {$\{i,j\}$};

\end{tikzpicture}
}
\caption{}
\end{figure}

Note that Definition \ref{def:gln foam coloring} implies that $F_i(c)$ and $F_{ij}(c)$ are closed orientable surfaces, and, in particular, they have even Euler characteristic. 

\begin{definition}
\label{def:rings}
Define the following rings. 
\begin{itemize}
 \item $\RN = \Z[x_1, \ldots, x_N]$. 
    
    \item $\RNsym \subset \RN$ is the subring of symmetric polynomials. 
    
    \item $\RN'=\RN[ (x_i - x_j)^{-1} \mid 1 \leq i < j \leq N]$ the extension of $\RN$ obtained by inverting $x_i-x_j$ for all $i<j$. We have an inclusion of rings $\RNsym \subset \RN \subset \RN'$.  These rings are all graded by setting $\deg(x_i) = 2$. 
\end{itemize}
\end{definition}

We are now ready to introduce Robert-Wagner closed foam evaluation \cite[Definition 2.12]{RWfoamev}. 

\begin{definition}
\label{def:gln foam evaluation} 
Let $F$ be a foam and $c\in \adm(F)$ a coloring. Define
\begin{align*}
    s(F,c) &= \sum_{i=1}^N i \chi(F_i(c))/2 + \sum_{1\leq i<j\leq N} \theta^+_{ij}(c), \\
    P(F,c) &= \prod_{f\in f(F)} P_f(c(f)), \\
    Q(F,c) &= \prod_{1\leq i < j \leq N} (x_i - x_j)^{\chi(F_{ij}(c))/2},\\
    \brak{F,c}_{\rm RW} &= (-1)^{s(F,c)} \frac{P(F,c)}{Q(F,c)}.
\end{align*}
Finally, the evaluation of $F$ is given by $\brak{F}_{\rm RW} = \sum\limits_{c\in \adm(F)} \brak{F,c}_{\rm RW}$
\end{definition}

Each $\brak{F,c}_{\rm RW}$ may have denominators, so, a priori, $\brak{F}_{\rm RW}$ is valued in $\RN'$. However, Robert-Wagner  prove that $\brak{F}_{\rm RW}$ lies in the subring   $\RNsym$ of symmetric polynomials \cite[Proposition 2.19]{RWfoamev}.

\section{Anchored \texorpdfstring{$\mathfrak{gl}_N$}{gln} link homology} 
\label{sec:anchored gln}
 Throughout this section, a positive integer $N$ will be fixed, and all link components are colored by an integer in $\{0, \ldots, N\}$. Recall the rings $\RNsym, \RN, \RN'$ from Definition \ref{def:rings}. We introduce the following additional notation.

 We introduce the following additional notation. 
 \begin{itemize}
 \item $\til{\RN} = \RN[ \sqrt{x_i - x_j} \mid 1\leq i < j \leq N ]$ the extension of $\RN$ obtained by introducing square roots of $x_i - x_j$ for $i<j$. 
    
    \item $\til{\RN}' = \til{\RN}[ (x_i - x_j)^{-1} \mid 1 \leq i < j \leq N] $ the ring obtained by inverting $x_i - x_j$, $i<j$, in $\til{\RN}$. 
 \end{itemize} 
 
 The rings $\til{\RN}, \til{\RN}'$ are graded by setting variables $x_i$ to be degree two.  Inclusions between these  rings are summarized in the diagram \eqref{eq:ring inclusion}. 
    
    \begin{equation}
    \label{eq:ring inclusion}
    \begin{tikzcd}[row sep = .3em, column sep =.2em ]
    & & \til{\RN} & \subset & \til{\RN}' \\
    & & \cup & & \cup \\
    \RNsym & \subset & \RN & \subset & \RN'
    \end{tikzcd}
    \end{equation}

\subsection{Anchored \texorpdfstring{$\mathfrak{gl}_N$}{gln} foams and their evaluations}
\label{sec:anchored gln foam evaluation}

In this section we extend the Robert-Wagner foam evaluation from Section \ref{sec:RW foam evaluation} to the annular setting.  Let $\L =   \{(0,0)\} \times \RR \subset \RR^3$ denote the  $z$-axis. We call $\L$ the \emph{anchor line}.

\begin{definition}
An \emph{anchored $\gl_N$ foam} is a $\gl_N$ foam $F$ such that intersections of $F$ with $\L$ occur transversely in the interior of facets of $F$. An intersection point of $F$ with $\L$ is called an \emph{anchor point}, and we let $\anch(F) = F\cap \L$ denote the set of anchor points. For $p\in \anch(F)$ lying on a facet $f$, we define its \emph{thickness} $\thick(p)$ to be the thickness of $f$.  

Anchored foams must also come equipped with a \emph{label} of each anchor point $p$, which consists of a subset $\l(p) \subset \NN = \{1,\ldots, N\}$ of cardinality equal to the thickness of the label on which $p$ lies.

The \emph{underlying foam} of $F$ is the $\gl_N$ foam obtained by forgetting anchor points and their labels. 

\end{definition}

A \emph{coloring} of an anchored foam means a coloring, in the sense of Definition \ref{def:gln foam coloring}, of the underlying foam. Let $F$ be an anchored foam and $c$ a coloring of $F$. For an anchor point $p\in \anch(F)$ lying on a facet $f$, its color $c(p)$ is defined to be the color of the facet, $c(p):=c(f)$.

We establish some notation before introducing anchored foam evaluation. 
\begin{itemize}
    \item For $A\subset \NN$, let $\b{A} = \NN\setminus A$ denote its complement. 
    \item For $A,B \subset \NN$, let $\less{A}{B}$ denote the subset of $A\times B$ consisting of pairs $(i,j)$ with $i\leq j$. Likewise, let $\lessstrict{A}{B}$ denote the subset of $A\times B$ consisting of pairs $(i,j)$ with $i<j$.
    \item For $A,B \subset \NN$, set 
    \[
      \OP(A,B) := \prod_{(i,j) \in \less{A}{B}}( x_i - x_j ) \prod_{(i,j) \in \less{B}{A}} (x_i - x_j).
      \]
\end{itemize}

\begin{lemma}
  For $A, B,C \subset \NN$, we have 
  \begin{enumerate}
    \item $\OP(A,B) =0 $ if and only if $A\cap B\neq\varnothing$. 
      \item $\OP(A,B) = (-1)^{\#\less{B}{A}} \prod_{(i,j) \in A\times B} (x_i - x_j)$.
      \item $\OP(A,B) = \OP(B,A)$.
      \item If $B\cap C=\varnothing$, then $\OP(A, B\cup C) = \OP(A,B)\OP(A,C).$ 
  \end{enumerate}
  \end{lemma}
  \begin{proof}
 All statements are immediate from the definition. 
  \end{proof}

\begin{definition}
Let $F$ be an anchored foam and $c$ a coloring of $F$. For $p\in \anch(F)$, define $\til{Q}(F,c,p) = \OP\left(c(p),\b{\ell(p)}\right)$, and set
\begin{align}
     \til{Q}(F,c) &= \left(\prod_{p\in \anch(F)} \til{Q}(F,c,p)\right)^{1/2}, \\
     \brak{F,c} &= (-1)^{s(F,c)}\frac{P(F,c)}{Q(F,c)} \cdot \til{Q}(F,c), \label{eq:brak F c}
\end{align}
where $s(F,c), P(F,c)$, and $Q(F,c)$ are as in Definition \ref{def:gln foam evaluation}. 
\end{definition}

Let us pause to comment on the above definition. As in Definition \ref{def:gln foam evaluation}, the evaluation $\brak{F}$ of an anchored foam $F$ will be  the sum of $\brak{F,c}$ over all colorings. Note that if $c(p)\neq \ell(p)$ for some $p\in \anch(F)$ then $\til{Q}(F,c,p) = \til{Q}(F,c) = \brak{F,c} = 0$. In light of this, we restrict the set of admissible colorings as follows.

\begin{definition}
For an anchored foam $F$, let $\adm(F)$ denote the set of colorings $c$ of $F$ such that $c(p) = \ell(p)$ for each $p\in \anch(F)$. Define
\[
\brak{F} = \sum_{c\in \adm(F)} \brak{F,c}. 
\]
\end{definition}

Due to the presence of the square root in $\til{Q}(F,c)$, a priori we have $\brak{F,c} \in \til{\RN}'$ (see the diagram \eqref{eq:ring inclusion} and the discussion above it for definitions of various rings). The following lemma shows that no square roots appear. 

\begin{lemma}
\label{lem:no square roots}
For an anchored foam $F$ and $c\in \adm(F)$, we have $\brak{F,c} \in \RN'$.
\end{lemma}
\begin{proof}
For $1\leq i< j\leq N$, a factor of $x_i-x_j$ appears under the square root in the definition of $\til{Q}(F,c)$ when either $i\in \l(p) = c(p), j\not\in \l(p)= c(p)$ or $j\in \l(p) = c(p), i\not\in \l(p)= c(p)$. Thus the power of $x_i - x_j$ appearing under the square root is equal to the number of intersection points between $F_{ij}(c)$ and $\L$, which is even since $F_{ij}(c)$ is a closed surface in $\RR^3$. 
\end{proof}

\begin{remark}
\label{rem:til Q gln}
In view of the proof of Lemma \ref{lem:no square roots}, we may write 
\[
\til{Q}(F,c) = \prod_{1\leq i<j\leq N} (x_i - x_j)^{\#(F_{ij}(c) \cap \L)/2}.  
\]
Letting $\anch(i,j)$ be the number of anchor points of $F$ which contain exactly one of $i$ or $j$ in their labels, we have $\#(F_{ij}(c) \cap \L) = \anch(i,j)$, so that the above expression of $\til{Q}(F,c)$ is independent of the admissible coloring $c$.  
\end{remark}

The remainder of this subsection is devoted to proving Proposition \ref{prop:no denominators}, which states that no denominators appear in $\brak{F}$. Let $S_N$ denote the symmetric group on $N$ letters. Fix an anchored foam $F$ and a permutation $\sigma \in S_N$. Let $\sigma(F)$ be the anchored foam whose underlying foam is the same as $F$ but the label of each anchor point $p$ is permuted by $\sigma$: if $\l_F(p)$ denotes the label of $p$ in $F$, then the label $\l_{\sigma(F)}(p)$ of $p$ in $\sigma(F)$ is given by $\l_{\sigma(F)}(p) = \sigma(\l(p))$. Likewise, if $c\in \adm(F)$ is admissible, then $\sigma(c) \in \adm(\sigma(F))$ denotes the coloring which colors a facet $f$ by $\sigma(c(f))$. Note that $S_N$ acts on the rings in diagram \eqref{eq:ring inclusion} by permuting the indices of variables. 

Define 
\[
\eps(F,\sigma) = \sum_{\substack{1\leq i < j \leq N, \\ \sigma(i) > \sigma(j)}} \anch(i,j)/2.
\]

\begin{lemma}
\label{lem:symmetric group action}
With the notation established above, we have
\[
\sigma\left(\brak{F,c}\right) = (-1)^{\eps(F,\sigma)}  \brak{\sigma(F), \sigma(c)}.
\]
In particular, $\brak{\sigma(F)} = \pm \brak{F}$, where the sign depends only on labels of anchor points of $F$ and on $\sigma$. 
\end{lemma}

\begin{proof}
By \cite[Lemma 2.17]{RWfoamev}, we have
\[
\sigma\left((-1)^{s(F,c)} \frac{P(F,c)}{Q(F,c))} \right)= (-1)^{s(\sigma(F),\sigma(c))} \frac{P(\sigma(F),\sigma(c))}{Q(\sigma(F),\sigma(c))}. 
\]
Using the reformulation of $\til{Q}(F,c)$ in Remark \ref{rem:til Q gln}, we see that $\sigma(\til{Q}(F,c)) = (-1)^{\eps(F,\sigma)} \til{Q}(\sigma(F),\sigma(c))$, which completes the proof.   
\end{proof}

For $1\leq i < j \leq N$, consider the subring 
\[
R_{i,j} = \RN [ (x_k - x_\l)^{-1} \mid 1 \leq k < \l \leq N, (k,\l)\neq (i,j) ]
\]
of $\RN'$. We also set $R_{j,i}=R_{i,j}$. We have
\[
\bigcap_{1\leq i < j \leq N} R_{i,j} = \RN,
\]
and $\sigma\in S_N$ sends $R_{i,j}$ isomorphically onto $R_{\sigma(i), \sigma(j)}$. 

We recall the following local operation on colors from \cite{RWfoamev}, called a \emph{Kempe move}. Let $F$ be an anchored foam,  $c\in \adm(F)$ a coloring, and $\Sigma \subset F_{ij}(c)$ a (not necessarily connected) sub-surface which is disjoint from $\L$. Define a coloring $c'\in \adm(F)$ to agree with $c$ on those  facets which are not in $\Sigma$, and for facets in $\Sigma$ the coloring $c'$ differs from $c$ by exchanging  $i$ and $j$. We say $c'$ is obtained from $c$ by an $(i,j)$-Kempe move.

\begin{proposition}
\label{prop:no denominators}
For any anchored foam $F$, we have $\brak{F} \in \RN$.
\end{proposition}

\begin{proof}
The argument is similar to that of \cite[Proposition 2.19]{RWfoamev}. First note that it suffices to show $\brak{F}\in R_{1,2}$. To see this, let $i<j$ and take $\sigma \in S_N$ with $\sigma(1) = i$, $\sigma(2) = j$. If $\brak{F}\in R_{1,2}$ for any anchored foam $F$, then by Lemma \ref{lem:symmetric group action} we have
\[
\sigma \left(\brak{\sigma^{-1}(F)} \right) = \pm \brak{F} \in R_{i,j},
\]
which implies $\brak{F}\in \RN$. 

Let us now show that $\brak{F} \in R_{1,2}$. Partition $\adm(F)$ as follows: the equivalence class $C_c$ of $c\in \adm(F)$ is the set of colorings obtained by $(1,2)$-Kempe moves on components of $F_{12}(c)$ which are disjoint from $\L$. We will show that 
\[
\sum_{c'\in C_c} \brak{F,c} \in R_{1,2},
\]
which will complete the proof. 

Write $F_{12}(c) = \Sigma_{\anch} \sqcup \Sigma_1 \sqcup \cdots \sqcup \Sigma_r$, where $\Sigma_1, \ldots, \Sigma_r$ are connected and disjoint from $\L$, and each component of $\Sigma_{\anch}$ intersects $\L$. Let $\sigma\in S_N$ denote the transposition $(1\ 2)$. Introduce the following notation
\begin{align*}
    P_{s}(F,c) &= \prod_{f \text{ a facet in } \Sigma_s} P(c(f)) \\
    P'(F,c) &= \prod_{f \text{ not a facet in } \cup_{s=1}^r \Sigma_s} P(c(f)) \\
    \h{Q}(F,c) &= \frac{Q(F,c) \prod\limits_{k\geq 3, s} (x_1 - x_k)^{\l_{\Sigma_s}(c,k)/2}}{(x_1 - x_2)^{\chi(F_{12}(c))/2}} \\
    \h{P}_{s}(F,c) &= P_s(F,c) \prod_{k\geq 3} (x_1 - x_k)^{\l_{\Sigma_s}(c,k)/2}\\
    T_s(F,c) &= \h{P}_s(F,c) + (-1)^{\chi(\Sigma_s)/2} \sigma \left( \h{P}_s(F,c)\right),
\end{align*}
where $\l_{\Sigma_s}(c,k)$ is the (even) integer defined in \cite[Lemma 2.10]{RWfoamev}. 

Arguing as in the proof of \cite[Proposition 2.19]{RWfoamev}, we have 
\[
\sum_{c'\in C_c} \brak{F,c}  = (-1)^{s(F,c)} \frac{P'(F,c)}{\h{Q}(F,c)}  \left( \prod_{s=1}^r (x_1 - x_2)^{-\chi(\Sigma_s)/2} T_s(F,c) \right) \frac{\til{Q}(F,c)}{(x_1 - x_2)^{\chi(\Sigma_{\anch})/2} }. 
\]
The term $\h{Q}(F,c)$ is not divisible by $x_1 - x_2$. Therefore denominators of the form $x_1 - x_2$ can appear in one of the two following ways. The first is when some $\Sigma_s$ is a $2$-sphere. In this case $T_s(F,c)$ is antisymmetric in $x_1, x_2$ and hence divisible by $x_1-x_2$, allowing to cancel this factor of $x_1-x_2$ in the denominator. The second is when some component $\Sigma'\subset \Sigma_{\anch}$ is a $2$-sphere. In this case $\Sigma'$ contains at least two anchor points $p_1, p_2$, each containing either $1$ or $2$ (but not both) in their labels. Thus their contribution 
\[
\sqrt{\til{Q}(F,c,p_1) \til{Q}(F,c,p_2)}
\]
to $\til{Q}(F,c)$ is divisible by $x_1 - x_2$, allowing to cancel with this contribution of $x_1-x_2$ to the denominator. 
\end{proof}

Closed anchored foams carry two types of gradings, which will induce gradings on state spaces.

\begin{figure}
\centering    
\subcaptionbox{$s(p) = 1$ \label{fig:positive intersection gln}}[.45\linewidth]
{\begin{tikzpicture}

\draw (0,0) -- (2.5,0);
\draw (0,0) -- (1,1);
\draw (1,1) -- (3.5,1);
\draw (2.5,0) -- (3.5,1);

\draw[red,preaction={draw, line width=3pt, white}] (1.75,1.5) -- (1.75,.5);
\draw[red, densely dashed] (1.75,.5) -- (1.75,0);
\draw[red] (1.75,0) -- (1.75,-.5);

\node(p) at (1.75,.5) {$*$};
\node[right] at (p) {$A$}; 

\draw[red,-{Stealth[scale=1.5]}] (1.75,1.3) -- (1.75,1.4);

\draw[->] (.9,.5) arc (110:430: .2 and .15);

\end{tikzpicture} 
}
\subcaptionbox{$s(p)=-1$ \label{fig:negative intersection gln}}[.45\linewidth]
{\begin{tikzpicture}

\draw (0,0) -- (2.5,0);
\draw (0,0) -- (1,1);
\draw (1,1) -- (3.5,1);
\draw (2.5,0) -- (3.5,1);

\draw[red,preaction={draw, line width=3pt, white}] (1.75,1.5) -- (1.75,.5);
\draw[red, densely dashed] (1.75,.5) -- (1.75,0);
\draw[red] (1.75,0) -- (1.75,-.5);

\node(p) at (1.75,.5) {$*$};
\node[right] at (p) {$A$}; 

\draw[red,-{Stealth[scale=1.5]}] (1.75,1.3) -- (1.75,1.4);

\draw[<-] (.9,.5) arc (110:430: .2 and .15);

\end{tikzpicture}
}
\caption{The oriented intersection number between an anchor point and $\L$.}
\label{fig:oriented intersection gln}
\end{figure}
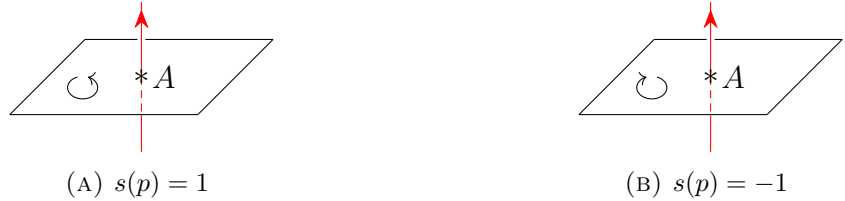

\begin{definition}
\label{def:qdeg and adeg gln}
The \emph{quantum grading}, denoted $\qdeg$, of an anchored foam $F$ is defined as 
\begin{equation*}
    \label{eq:qdeg}
    \qdeg(F) = \deg(F^{\rm un}) + \sum_{p \in \anch(F)} \thick(p)(N-\thick(p)) ,
\end{equation*}
where $F^{\rm un}$ is the underlying foam and $\deg(F^{\rm un})$ is as in Definition \ref{def:qdeg of gln foam}. 

Anchored foams carry an additional $\Z^N$-grading, called the \emph{annular degree}. Let $w_1, \ldots, w_N$ denote the standard basis of $\Z^N$. Given $A\subset \NN$, set 
\[
w_A = \sum_{i\in A} w_i.
\]
For
an anchored foam $F$, we define $\adeg(F) \in\Z^N$ as follows. Orient the anchor line $\L$ from bottom to top. For an anchor point $p\in \anch(F)$ lying on a facet $f$, denote by $s(p) \in \{\pm 1\}$ the oriented intersection number of $f$ and $\L$ at $p$ (see Figure \ref{fig:oriented intersection gln}). Set 
\begin{equation*}
\label{eq:annular degree}
    \adeg(F) = \sum_{p\in \anch(F)} s(p) w_{\ell(p)}. 
\end{equation*}
\end{definition}

We denote the usual $\Z$-grading on the ground ring $\RN$ by $\qdeg$, given by $\qdeg(x_i) = 2$. Define the trivial\footnote{While the ground ring is trivially graded with respect to $\adeg$, state spaces of annular webs will, in general, have non-trivial annular gradings.} $\Z^N$-grading on $\RN$ by setting $\adeg(x_i) = (0,\ldots, 0) $. We view $\RN$ as $\Z\oplus \Z^N$-graded via $(\qdeg, \adeg)$. The following lemma says that anchored $\gl_N$ foam evaluation respects degrees. 

\begin{lemma}
\label{lem:degrees of closed anchored gln foams}
If $F$ is a closed anchored foam with an admissible coloring, then $\qdeg(F) = \qdeg(\brak{F})$ and $\adeg(F) = \adeg(\brak{F}) = 0$. 
\end{lemma}

\begin{proof}
Let $c \in \adm(F)$. The first statement is clear from the formula for $\brak{F,c}$, equation \eqref{eq:brak F c}.  The argument for the second statement is similar to the proof of \cite[Proposition 4.23]{AkhKh}. Consider the intersection of $F$ with a generic half-plane with boundary $\L$, resulting in a web $\Gamma$ with boundary on $\L$. Each edge $e$ in $\Gamma$  is the intersection of the half-plane with a facet $f$ of $F$, and we set $c(e):=c(f)$. Edges with boundary on $\L$ are colored according to their label. For each internal vertex $v$ of $\Gamma$, with edges $e_1, e_2, e_3$ incident to $v$, let
\[
w(v) := \eps_1 \sum_{i\in c(e_1) } w_i + \eps_2 \sum_{i\in c(e_2) } w_i + \eps_3 \sum_{i\in c(e_3) } w_i,
\]
where $\eps_j = 1$ if $e_j$ is oriented towards $v$ and $\eps_j =-1$ if $e_j$ is oriented away from $v$.  Each $w(v) = 0$ since $c$ is a coloring of $F$, so 
\begin{equation}
\label{eq:vertex contributions}
\sum_{v \text{ a vertex of } \Gamma}  w(v) = 0.
\end{equation}
Finally, observe that $\adeg(F)$ is equal to the above sum. To see this, let $E$ be the set of edges of $\Gamma$ with exactly one endpoint on $\L$. For $e\in E$, let $w(e) = \eps \sum_{i\in c(e)} w_i$, where $\eps = -1$ if $e$ is oriented towards $\L$ and $+1$ otherwise. The sum in \eqref{eq:vertex contributions} is equal to $\sum_{e \in E} w(e)$, since each internal edge contributes the same quantity with opposite sign to \eqref{eq:vertex contributions}, leaving only contributions from edges in $E$. On the other hand, we can compute $\adeg(F)$ from labels and orientations on boundary points of $\Gamma$. Two boundary points which are connected by an edge contribute the same quantity to $\adeg(F)$ with opposite signs, and the remaining boundary points contribute precisely  $\sum_{e \in E} w(e)=0$.

\end{proof}

\subsection{Local relations} 
\label{sec:anchored gln local relations}
This section establishes local relations satisfied by anchored $\gl_N$ foam evaluation.

\begin{lemma}
The local relation \eqref{eq:polynomial removal} holds. 
\begin{equation}
\label{eq:polynomial removal}
    \begin{aligned}
    \begin{tikzpicture}[scale=.9]

\draw (0,0) -- (2.5,0);
\draw (0,0) -- (1,1);
\draw (1,1) -- (3.5,1);
\draw (2.5,0) -- (3.5,1);

\draw[red,preaction={draw, line width=3pt, white}] (1.75,1.5) -- (1.75,.5);
\draw[red, densely dashed] (1.75,.5) -- (1.75,0);
\draw[red] (1.75,0) -- (1.75,-.5);

\draw[blue, ->] (-.2,.5) -- (1,.5); 

\node[left] at (-.2,.5) {$P$};

\node(p) at (1.75,.5) {$*$};
\node[anchor= west] at (p) {$A$};

\begin{scope}[shift = {(5.5,0)}]

\draw (0,0) -- (2.5,0);
\draw (0,0) -- (1,1);
\draw (1,1) -- (3.5,1);
\draw (2.5,0) -- (3.5,1);

\draw[red,preaction={draw, line width=3pt, white}] (1.75,1.5) -- (1.75,.5);
\draw[red, densely dashed] (1.75,.5) -- (1.75,0);
\draw[red] (1.75,0) -- (1.75,-.5);

\node(p) at (1.75,.5) {$*$};
\node[anchor= west] at (p) {$A$};
\end{scope}

\node at (4.5,.5) {$=\; P(A)$};

\end{tikzpicture}
    \end{aligned}
\end{equation}
\end{lemma}
\begin{proof}
This is immediate from the definition. 
\end{proof}

\begin{lemma}
\label{lem:theta foam gln}
Consider the foam $F$ shown in \eqref{eq:theta foam anchored general}.
\begin{equation}
    \label{eq:theta foam anchored general}
    \begin{aligned}
    \begin{tikzpicture}[scale=1.5,decoration={
    markings,
    mark=at position 0.75 with {\arrow{Stealth[scale=2]}}}]

\draw (0,0) circle [radius=1];

\draw[densely dashed,fill=gray!40] (0,0) ellipse (1 and .3);
\draw[densely dashed] (1,0) arc(0:180:1 and .3);
\draw[postaction={decorate}] (1,0) arc(0:-180:1 and .3);

\draw[red] (0,1.5) -- (0,1);
\draw[densely dashed,red] (0,1) -- (0,-1);
\draw[red] (0,-1) -- (0,-1.5);

\node(dt) at (.4,.6) {$b$};
\node(dm) at (.4,0) {$a+b$};
\node(db) at (.4,-.6) {$a$};
j

\node(t) at (0,1) {$*$};
\node(m) at (0,0) {$*$};
\node(b) at (0,-1) {$*$};

\node[anchor = south east] at (t) {$B$};
\node[anchor = east] at (m) {$C$}; 
\node[anchor = north east] at (b) {$A$}; 

\node(lt) at (-1.5,.6) {$P_2$};
\node(lm) at (-1.5,0) {$P_3$};
\node(lb) at (-1.5,-.6) {$P_1$};

\draw[thick, ->, blue] (lt) -- (-.6,.6);
\draw[thick, ->, blue] (lm) -- (-.6,0);
\draw[thick, ->, blue] (lb) -- (-.6,-.6);
\end{tikzpicture}
    \end{aligned}
\end{equation}
Then 
\[
\brak{F} = 
\begin{cases}
(-1)^{\sum\limits_{i\in C} i + \#\lessstrict{B}{A}  } \ P_1(A) P_2(B) P_3(C) & \text{ if } A\cup B=C, \\ 
0 & \text{ otherwise}.
\end{cases}
\]
\end{lemma}

\begin{proof}
If $C\neq A\cup B$ then $F$ has no admissible colorings. Suppose then that $C = A\cup B$, in which case $A$ and $B$ are disjoint, and there is a single admissible coloring $c$. We have $P(F,c) = P_1(A) P_2(B) P_3(C)$. For $i<j$, the bicolored surface $F_{ij}(c)$ is either empty or a $2$-sphere. The latter occurs in two cases: first, when $(i,j)$ or $(j,i)$ is in $A \times B$, where the $2$-sphere is the union of the thickness $a$ and $b$ facets. The second case when $(i,j)$ or $(j,i)$ is in $C \times \b{C}$, where the $2$-sphere is the union of the thickness $a+b$ facet and either the thickness $a$ or the thickness $b$ facet. Thus 
\[
Q(F,c) = \OP(A,B) \OP(C,\b{C}).  
\]
Contributions from anchor points also equal the above term:
\begin{align*}
\til{Q}(F,c) & = \left(\OP(A,\b{A}) \cdot \OP(B,\b{B}) \cdot \OP(C,\b{C}) \right)^{1/2} \\    
            &= \left(\OP(A,\b{C}) \cdot \OP(A,B) \cdot \OP(B, \b{C}) \cdot \OP(B,A) \cdot \OP(C,\b{C}) \right)^{1/2} \\
            &= \OP(A,B) \OP(C,\b{C}). 
\end{align*}

It remains to compute the sign. For $i\in \NN$, the monochrome surface $F_i(c)$ is empty whenever $i \not\in C$, and otherwise for $i\in C$, $F_i(c)$ is a $2$-sphere. Letting $Z$ denote the depicted singular circle and fixing $i<j$, we see that $Z$ is positive with respect to $(i,j)$ if and only if $i\in B$, $j\in A$. So we obtain 
\[
s(F,c) = \sum_{i\in C} i + \#\lessstrict{B}{A}. 
\]
\end{proof}

The following lemmas will be crucial for identifying state spaces assigned to annular webs. 

\begin{lemma}
\label{lem:gln sphere}
Consider the foam $F$ shown in \eqref{eq:sphere anchored}, consisting of a thickness $a$ sphere decorated by the symmetric polynomial $P$ and which intersects $\L$ twice, with anchor points labeled $A$ and $B$. 
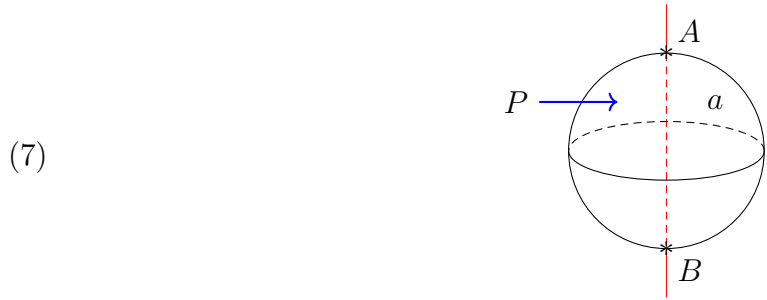
\begin{equation}
    \label{eq:sphere anchored}
    \begin{aligned}
    \begin{tikzpicture}[scale=1.3]
\draw[red] (0,1.5) -- (0,1);
\draw[densely dashed,red] (0,1) -- (0,-1);
\draw[red] (0,-1) -- (0,-1.5);

\draw (0,0) circle [radius=1];


\draw[densely dashed] (1,0) arc(0:180:1 and .3);
\draw (1,0) arc(0:-180:1 and .3);

\draw[thick, ->, blue] (-1.3,.5) -- (-.5,.5);
\node[left] at (-1.3,.5) {$P$};

\node at (.5,.5) {$a$}; 

\node(t) at (0,1) {$*$};
\node(b) at (0,-1) {$*$};

\node[anchor = south west] at (t) {$A$}; 
\node[anchor = north west] at (b) {$B$}; 

\end{tikzpicture}
    \end{aligned}
\end{equation}
Then 
\[
\brak{F} = 
\begin{cases}
(-1)^{\sum\limits_{i \in A} i} P(A) & \text{ if } A = B, \\ 
0 & \text{ otherwise}.
\end{cases}
\]
\end{lemma}

\begin{proof}
If $A\neq B$ then there are no admissible colorings, and if $A=B$ there is exactly one admissible coloring $c$ which colors $F$ by $A$. We have $s(F,c) = \sum_{i\in A} i$, $P(F,c) = P(A)$, and $Q(F,c) = \OP(A,\b{A}) = \til{Q}(F,c)$. 
\end{proof}

\begin{lemma}
\label{lem:neck cutting gln}
The local relation \eqref{eq:essential circle id relation} holds.
\begin{equation}
    \label{eq:essential circle id relation}
    \begin{aligned}
    \begin{tikzpicture}[decoration={
    markings,
    mark=at position 0.7 with {\arrow{Stealth[scale=1]}}}]


\draw[postaction={decorate}] (2,0) arc (0:-180:1 and .3);
\draw[densely dashed] (2,0) arc (0:180:1 and .3);

\draw[postaction={decorate}] (2,4) arc (0:-180:1 and .3);
\draw[] (2,4) arc (0:180:1 and .3);

\draw (0,0) -- (0, 4);
\draw (2,0) -- (2, 4);

\node at (1.5,2) {$a$};


\draw[red,preaction={draw, line width=3pt, white}] (1,4.5) --++ (0,-.75);
\draw[red,densely dashed] (1,3.65) -- (1,-.28);
\draw[red] (1,-.35) --++(0,-.3);


\begin{scope}[shift={(5.75,0)}]


\draw[postaction={decorate}] (2,0) arc (0:-180:1 and .3);
\draw[densely dashed] (2,0) arc (0:180:1 and .3);

\draw (2,0) arc (0:180:1 and 1.2);


\draw[postaction={decorate}] (2,4) arc (0:-180:1 and .3);
\draw[] (2,4) arc (0:180:1 and .3);

\draw (2,4) arc (0:-180:1 and 1);

\draw[red,preaction={draw, line width=3pt, white}] (1,4.5) --++ (0,-.75);;
\draw[red,densely dashed] (1,3.65) -- (1,3);
\draw [red] (1,3) -- (1,2.45);
\draw[red] (1,2.45) -- (1,1.25);
\draw[red,densely dashed] (1,1.25) --(1,-.28);
\draw[red] (1,-.35) --++(0,-.3);


\node (a1) at (1,3) {$*$};

\node (a2) at (1,1.25) {$*$};

\node[anchor = north west] at (a1) {$\mathsmaller{A}$};
\node[anchor=south west] at (a2) {$\mathsmaller{A}$};

\node at (1.5,3.5) {$\mathsmaller{a}$};
\node at (1.5,.5) {$\mathsmaller{a}$};

\end{scope}

\node at (4,2) {$\displaystyle =\ \  {\sum_{\mathsmaller{\substack{A \subset [N] \\ \# A = a }} }} (-1)^{\sum\limits_{i\in A} i} $};
\end{tikzpicture}
    \end{aligned}
\end{equation}
\end{lemma}

\begin{proof}
Denote by $F$ the foam on the left-hand side on the equality. For each $A\subset \NN$ of order $a$, denote by $\adm_A(F)$ the admissible colorings of $F$ in which the depicted annulus is colored by $A$, and let $G^A$ be the foam summand on the right-hand side of the equality whose top anchor point is labeled $A$. There is a natural bijection $\adm_A(F) \cong \adm(G^A)$. For $c\in \adm_A(F)$ with corresponding coloring $c'\in \adm(G^A)$, we will show that 
\[
\brak{F,c} = (-1)^{\sum_{i\in A} i} \brak{G^A,c'},
\]
from which relation \eqref{eq:essential circle id relation} follows. 

For $1\leq i \leq N$, we have 
\[
\chi(G^A_i(c')) =
\begin{cases}
\chi(F_i(c)) +2 & \text{ if } i \in A, \\
\chi(F_i(c)) & \text{ if } i \not\in A.
\end{cases}
\]
It follows that $s(F,c) = s(G^A, c') + \sum_{i\in A} i$.

Next, we have $P(G^A,c') = P(F,c)$. For $i<j$, we have 
\[
\chi(G^A_{ij}(c')) = 
\begin{cases}
\chi(F_{ij}(c)) + 2 & \text{ if exactly one of } i,j \text{ is in } A, \\
\chi(F_{ij}(c)) & \text{ otherwise}.

\end{cases}
\]
so 
\[
Q(G^A,c') = Q(F,c)\cdot \OP(A,\b{A}). 
\]
On the other hand, we have 
\[
\til{Q}(G^A,c')= \til{Q}(F,c) \cdot \OP(A, \b{A}),
\]
which verifies $\brak{F,c} = (-1)^{\sum_{i\in A} i} \brak{G^A,c'}$ and completes the proof of the lemma. 

\end{proof}

\subsection{State spaces of annular \texorpdfstring{$\mathfrak{gl}_N$}{gln} webs} 
\label{sec:state spaces annular gln}

In this section we define state spaces for annular $\gl_N$ webs using universal construction applied to anchored foam evaluation. The main result is Theorem \ref{thm:anchored gln state space}, which identifies state spaces.

\begin{definition}
An \emph{annular $\gl_N$ web} is a $\gl_N$ web embedded in the punctured plane $\P$. Let $\Gamma_0, \Gamma_1$ be annular webs. An \emph{anchored cobordism} from $\Gamma_0$ to $\Gamma_1$ is a $\gl_N$ foam with boundary $F\subset \RR^2\times [0,1]$ (see Definition \ref{def:gln foam with boundary}) from $\Gamma_0$ to $\Gamma_1$ (viewing $\Gamma_0$ and $\Gamma_1$ as subsets of the plane $\RR^2$) such that $F$ intersects the line segment $\L_{[0,1]} := \{(0,0)\} \times [0,1]$ only in the interior of its facets. These intersection points, called anchor points as for closed foams, are required to be transverse, and each anchor point $p$ is labeled by a subset $\ell(p)\subset \NN$ of cardinality equal to the thickness of the facet on which $p$ lies. Write $\anch(F) = F\cap \L$ to denote the anchor points. 

Anchored cobordisms are considered up to ambient isotopy of $\RR^2 \times [0,1]$ which is the identity near $\d (\RR^2\times [0,1])$ and which maps $\L_{[0,1]}$ to itself. Define $\qdeg(F)$ and $\adeg(F)$ for an anchored cobordism $F$ as in Definition \ref{def:qdeg and adeg gln}. 

Let $\aFoam_N$ denote the category of annular webs and anchored cobordisms between them. The gradings $\qdeg$ and $\adeg$ are additive under composition of anchored cobordisms. 
\end{definition}

We define the state space of an annular web $\Gamma$ via universal construction in the usual way. Let $\Fr(\Gamma)$ denote the free $\RN$-module generated by all anchored cobordisms from the empty web to $\Gamma$. Consider the symmetric bilinear form 
\[
\brak{-,-} : \Fr(\Gamma) \times \Fr(\Gamma) \to \RN 
\]
given by $\brak{F,G} = \brak{\b{F}G}$, where $\b{F}$ is the reflection of $F$ through $\RR^2$. 
Define the \emph{state space} $\brak{\Gamma}$ to be the quotient of $\Fr(\Gamma)$ by the kernel 
\[
\ker(\brak{-,-}) = \{x\in \Fr(\Gamma) \mid \brak{x,y} = 0 \text{ for all } y\in \Fr(\Gamma)\}
\]
of this bilinear form.

Lemma \ref{lem:degrees of closed anchored gln foams} implies that $\qdeg$ and $\adeg$ descend to gradings on $\brak{\Gamma}$, so $\brak{\Gamma}$ is $\Z\oplus \Z^N$-graded. It follows immediately from the definitions that an anchored cobordism $F : \Gamma_0 \to \Gamma_1$ induces a map $\brak{F}: \brak{\Gamma_0} \to \brak{\Gamma_1}$  of degree $(\qdeg(F), \adeg(F))$, defined by sending the equivalence class of a basis element $V\in \Fr(\Gamma_0)$ to the equivalence class of  $F \circ V \in \Fr(\Gamma_1)$. Thus we obtain a functor 
\[
\brak{-} : \aFoam_N \to \RN \gNmod,
\]
where $\RN\gNmod$ denotes the category of $\Z\oplus\Z^N$-graded $\RN$-modules.

For a $\Z\oplus\Z^N$-graded $\RN$-module $M$ and $i\in \Z$, $j=(j_1, \ldots, j_N) \in \Z^N$, we will write  $M\{(i,j)\}$ to denote the module $M$ with gradings shifted up by $(i,j)$. Given $f\in \Z_{\geq 0} [q^{\pm 1}, a_1^{\pm 1}, \ldots, a_N^{\pm 1}]$, write $f= \sum\limits_{I = (i, j_1, \ldots, j_N) \in \Z^{N+1}} m_I q^{i} a_1^{j_1} \cdots a_N^{j_N}$, and set 
\[
M\{f\} := \bigoplus_{I = (i, j_1, \ldots, j_N) \in \Z^{N+1}} M^{\oplus m_I} \{(i, {j_1}, \ldots, {j_N})\}.
\]
A free graded $R_N$-module is of the form $\bigoplus_{i=1}^k R_N\{f_i\}$, and its graded rank is defined to be $f_1 + \cdots + f_k$.

We establish the following sequence of lemmas before proving Theorem \ref{thm:anchored gln state space}. 

\begin{lemma}
\label{lem:zero edge removal}
Let $\Gamma$ be an annular web, and let $\Gamma'$ be the annular web obtained from $\Gamma$ by deleting all the edges of thickness zero. Then there is an isomorphism of state spaces $\brak{\Gamma} \cong \brak{\Gamma'}$. 
\end{lemma}
\begin{proof}
The argument  in \cite[Claim 3.33]{RWfoamev} applies without modification. 
\end{proof}

\begin{definition}[{\cite[Notation 3.6, Definition 3.7]{RWfoamev}}]
Let $C$ be an oriented cycle in an annular web $\Gamma$, and let $\Gamma'$ be the annular web obtained from $\Gamma$ by reversing the orientation of each edge in $C$ and replacing its thickness $a$ by $N-a$. The web $\Gamma'$ is said to be obtained by a \emph{cycle move along $C$}. The cycle $C$ is \emph{face-like} if it bounds a disk in $\RR^2\setminus \Gamma$. 
\end{definition}

Note that face-like cycles are allowed to bound a disk which contains the puncture.

\begin{lemma}
\label{lem:cycle move}
If $\Gamma'$ is obtained from $\Gamma$ by a cycle move, then $\brak{\Gamma'} \cong \brak{\Gamma}$.
\end{lemma}

\begin{proof}
By \cite[Lemma 3.9]{RWfoamev}, it suffices to prove the claim when $C$ is face-like. If $C$ bounds a disk which does not contain the puncture, then the isomorphism can be taken to be the same one as in \cite[Claim 3.34]{RWfoamev}. Suppose that $C$ bounds a disk which contains the puncture. We claim that the relation \eqref{eq:annular face like} holds, where:
\begin{itemize}
    \item the shaded facets have thickness $N$,
    \item the bold dashed arrows are not seam lines indicating a meeting of facets, but rather the orientation on the edges of the web that would be created if the foam were sliced horizontally, 
    \item the relation holds not only for squares but for any face-like cycle, which is why most orientations and thicknesses are omitted.
\end{itemize}

\begin{equation}
\label{eq:annular face like} 
    \begin{aligned}
    \includestandalone{annular_face_like_cycle}
    \end{aligned}
\end{equation}

Since anchor points on thickness $N$ facets do not contribute to the evaluation, the identity \eqref{eq:annular face like}  follows from the proof of \cite[Lemma 3.21]{RWfoamev}. Therefore a modification of the foams witnessing the isomorphism in \cite[Claim 3.34]{RWfoamev}, given by passing the anchor line through the thicknesses $N$ facets, yields the desired isomorphism in the case where $C$ bounds a disk containing the puncture.

\end{proof}

Recall that $w_1, \ldots, w_N$ denote the standard basis vectors of $\Z^N$, and that if $A\subset \NN$, we set $w_A = \sum_{i\in A} w_i$.

\begin{theorem}
\label{thm:essential circle gln}
Let $\Gamma$ be an annular web containing a non-contractible counterclockwise oriented circle $Z$ of thickness $a$ which bounds a disk in $\RR^2\setminus \Gamma$. Let $\Gamma' = \Gamma\setminus Z$ denote the annular web obtained by removing $Z$ from $\Gamma$. There is an isomorphism $\brak{\Gamma} \cong  \bigoplus\limits_{\substack{A \subset \NN \\ \#A= a }} \brak{\Gamma'} \{(0,w_A)\}$. This is depicted in \eqref{eq:gln essential circle}. 
\begin{equation}
\label{eq:gln essential circle}
   \left\langle
\!\begin{gathered}
\begin{tikzpicture}[scale=.6,x=10mm,y=10mm,baseline={([yshift=-.5ex]current bounding box.center)},>=Stealth,decoration={
    markings,
    mark=at position 0.5 with {\arrow{Stealth}}}
    ] 
\node at (0,.75) {$\times$}; 
\draw[->] (0,0) -- (.1,0);
\draw[] (0,0) arc (-90:90:.75) node [midway,right] {$a$};
\draw (0,1.5) arc (90:270:.75);
\end{tikzpicture}
\end{gathered}\right\rangle\ 
\cong  \ \bigoplus\limits_{\substack{A \subset \NN \\ \#A= a }} \brak{\varnothing} \{(0,w_A)\}
\end{equation}
If the orientation of $Z$ is clockwise then the isomorphism holds with negated degree shifts. 
\end{theorem}

\begin{figure}
\centering
\subcaptionbox{The anchored cobordism $F^A : \Gamma \to \Gamma'$.
\label{fig:F^A}}[.45\linewidth]
{\begin{tikzpicture}[decoration={
    markings,
    mark=at position 0.7 with {\arrow{Stealth[scale=1]}}}]

\draw[postaction={decorate}] (0,0) arc (180:360:1 and .3);
\draw[densely dashed] (2,0) arc (0:180:1 and .3);

\draw (2,0) arc (0:180:1 and 1.2);

\draw [red] (1,2.5) -- (1,1.2);
\draw[red,densely dashed] (1,1.2) -- (1,0);
\node at (1,2.5) {$\color{red}\times$};

\node at (1,0) {$\color{red}\times$};


\node (a2) at (1,1.2) {$*$};


\node[anchor=south east] at (a2) {${A}$};


\node at (1.5,.6) {$\mathsmaller{a}$};

\end{tikzpicture}
}
\subcaptionbox{The anchored cobordism $G^A : \Gamma' \to \Gamma$.
\label{fig:G^A}}[.45\linewidth]
{\begin{tikzpicture}[decoration={
    markings,
    mark=at position 0.7 with {\arrow{Stealth[scale=1]}}}]

\draw[red] (1,4.5) --++ (0,-.25);
\draw[red,densely dashed] (1,4.15) -- (1,3.5);
\draw [red] (1,3.5) -- (1,2);

\node at (1,4.5) {$\color{red}\times$};

\node at (1,2) {$\color{red}\times$};
\draw[postaction={decorate},preaction={draw, line width=3pt, white}] (0,4.5) arc (180:360:1 and .3);
\draw[] (2,4.5) arc (0:180:1 and .3);

\draw (2,4.5) arc (0:-180:1 and 1.2);

\node (a1) at (1,3.3) {$*$};

\node[anchor = north east] at (a1) {${A}$};

\node at (1.5,3.9) {$\mathsmaller{a}$};

\end{tikzpicture}
}
\caption{The foams $F^A$ and $G^A$ in the proof of Theorem \ref{thm:essential circle gln}}\label{fig:essential circle foam maps gln}
\end{figure}
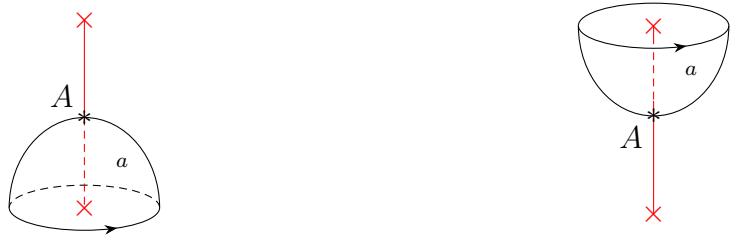

\begin{proof}
For $A \subset \NN$ with $\#A= a$, let $F^A : \Gamma \to \Gamma'$ denote the anchored cobordism depicted in Figure \ref{fig:F^A}, and let $G^A : \Gamma'\to \Gamma$ denote the anchored cobordism depicted in Figure \ref{fig:G^A} (only the relevant part of the foams are shown - outside of the depicted region they are the identity on $\Gamma'$). It is straightforward to verify that $\qdeg(F^A) = \qdeg(G^A) = 0$, $\adeg(F^A) = -w_A$, and $\adeg(G^A) = w_A$.

Define $\Phi: \brak{\Gamma} \to \brak{\Gamma'}$ to be the $\genfrac(){0pt}{1}{N}{a} \times 1$ matrix of foams with entries 
\[
(-1)^{\sum_{i\in A} i} \brak{F^A},
\]
and $\Psi: \brak{\Gamma'} \to \brak{\Gamma}$ to be the $1\times \genfrac(){0pt}{1}{N}{a}$ matrix of foams with entries $\brak{G^A}$. Lemma \ref{lem:neck cutting gln} and Lemma \ref{lem:gln sphere} imply, respectively, that $\Psi\Phi = \id_\Gamma$ and $\Phi\Psi = \id_{\Gamma'}$. 
\end{proof}

Consider any one of MOY relations in Figure \ref{fig:MOY relations}, where the depicted portions of webs are contained in disks disjoint from the puncture, and write it in the form $\Gamma = \sum_i q^{n_i} \Gamma_i$. By a \emph{categorified} MOY relation we mean an isomorphism of state spaces $\brak{\Gamma} \cong \bigoplus_i \brak{\Gamma_i} \{(n_i,0)\}$. 

\begin{lemma}
All the categorified MOY relations hold for anchored $\gl_N$ state spaces. 
\end{lemma}

\begin{proof}
Robert-Wagner show that there is an isomorphism of state spaces \[
\brak{\Gamma}_{\rm RW} = \bigoplus_i \brak{\Gamma_i}_{\rm RW} \{n_i\}.
\]
For most of the relations, these isomorphisms are realized explicitly via maps induced by foams: categorified versions of Figures \ref{MOYa}, \ref{MOYb}, \ref{MOYc}, \ref{MOYd}, and \ref{MOYg} follow from the relations \cite[Equation (10), Claim 3.35, Equation (12), Equation (13), Claim 3.36]{RWfoamev}, respectively. The remaining relations are Figures \ref{MOYe} and \ref{MOYf}. These too can be realized via foams. The relation in Figure \ref{MOYe} is a special case of Figure \ref{MOYg}. Figure \ref{MOYf} is a combination of specializing Figure \ref{MOYg} and applying a cycle move, and cycle moves can be realized by anchored foams which induce isomorphisms of state spaces by Lemma \ref{lem:cycle move}. The arguments in \cite{RWfoamev} establishing foam relations are completely local and apply unchanged in the annular setting when all the local foam relations are disjoint from $\L$. 
\end{proof}

\begin{theorem}
\label{thm:anchored gln state space}
Let $\Gamma$ be an annular web. Then $\brak{\Gamma}$ is a finitely generated, free, graded  $\RN$-module. Its graded rank can be computed by applying the MOY relations in Figure \ref{fig:MOY relations} (where the local pictures are all disjoint from the puncture) and the additional relations shown in \eqref{eq:annular MOY relation}. 

\begin{equation}
\label{eq:annular MOY relation}
    \begin{aligned}
    \begin{tikzpicture}[>=Stealth,decoration={
    markings,
    mark=at position 0.5 with {\arrow{Stealth}}}
    ] 
\begin{scope}[shift={(-6,0)}]
\node at (0,.5) {$\times$}; 

\draw[->] (0,0) -- (.1,0);

\draw[] (0,0) arc (-90:90:.5) node [midway,right] {$a$};
\draw (0,1) arc (90:270:.5);

\node at (2.5,.25) {$\  =  \ \displaystyle \sum_{\substack{A\subset [N]\\\#A = a }} \  \prod_{i\in A} a_i$};
\end{scope}


\node at (0,.5) {$\times$}; 

\draw[->] (0,0) -- (-.1,0);

\draw[] (0,0) arc (-90:90:.5) node [midway,right] {$a$};
\draw (0,1) arc (90:270:.5);

\node at (2.5,.25) {$\  = \ \displaystyle \sum_{\substack{A\subset [N] \\ \#A = a} } \  \prod_{i\in A} a_i^{-1}$};
\end{tikzpicture}
    \end{aligned}
\end{equation}

\end{theorem}

\begin{proof}
By Lemma \ref{lem:zero edge removal}, we may assume that $\Gamma$ has no edges of thickness zero. Suppose that all edges of $\Gamma$ are of thickness $1$ or $2$. In this case the underling graph of $\Gamma$, not including closed loops, is bipartite. By \cite[Lemma 3.17]{AkhKh}, $\Gamma$ contains either a closed innermost\footnote{Innermost means that the loop bounds a disk in $\RR^2\setminus \Gamma$ (the disk may contain the puncture).} loop, bigon, or square region. Every closed innermost loop can be removed by using either the categorified MOY relation Figure \ref{MOYa} or Theorem \ref{thm:essential circle gln}. There are two types of bigon faces in this case, which can be removed using the categorified MOY relations Figure \ref{MOYc} and Figure \ref{MOYd}. Finally, for this step, consider a square face of $\Gamma$. Up to rotation and reflection, there are three types of square faces: 
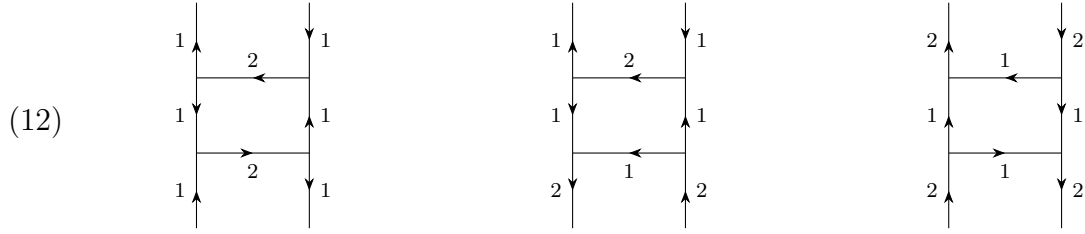
\begin{equation}
    \begin{aligned}
    \begin{tikzpicture}[decoration={
    markings,
    mark=at position 0.5 with {\arrow{Stealth}}}
    ]

\draw[postaction={decorate}] (0,0) --++ (0,1) node[midway, left] {$\mathsmaller{1}$}; 

\draw[postaction={decorate}] (1.5,0) --++ (-1.5,0) node[midway, above] {$\mathsmaller{2}$}; 

\draw[postaction={decorate}] (1.5,1) --++ (0,-1) node[midway, right] {$\mathsmaller{1}$}; 

\draw[postaction={decorate}] (0,0) --++ (0,-1) node[midway,left] {$\mathsmaller{1}$}; 

\draw[postaction={decorate}] (0,-1) --++ (1.5,0) node[midway,below] {$\mathsmaller{2}$};

\draw[postaction={decorate}] (1.5,-1) --++ (0,1) node[midway,right] {$\mathsmaller{1}$};

\draw[postaction={decorate}] (0,-2) --++ (0,1) node[midway,left] {$\mathsmaller{1}$}; 

\draw[postaction={decorate}] (1.5,-1) --++ (0,-1) node[midway,right] {$\mathsmaller{1}$};


\begin{scope}[shift={(5,0)}]

\draw[postaction={decorate}] (0,0) --++ (0,1) node[midway, left] {$\mathsmaller{1}$}; 

\draw[postaction={decorate}] (1.5,0) --++ (-1.5,0) node[midway, above] {$\mathsmaller{2}$}; 

\draw[postaction={decorate}] (1.5,1) --++ (0,-1) node[midway, right] {$\mathsmaller{1}$}; 

\draw[postaction={decorate}] (0,0) --++ (0,-1) node[midway,left] {$\mathsmaller{1}$}; 

\draw[postaction={decorate}] (1.5,-1) -- (0,-1) node[midway,below] {$\mathsmaller{1}$};

\draw[postaction={decorate}] (1.5,-1) --++ (0,1) node[midway,right] {$\mathsmaller{1}$};

\draw[postaction={decorate}] (0,-1) -- (0,-2) node[midway,left] {$\mathsmaller{2}$}; 

\draw[postaction={decorate}] (1.5,-2) -- (1.5,-1) node[midway,right] {$\mathsmaller{2}$};
\end{scope}


\begin{scope}[shift={(10,0)}]

\draw[postaction={decorate}] (0,0) --++ (0,1) node[midway, left] {$\mathsmaller{2}$}; 

\draw[postaction={decorate}] (1.5,0) --++ (-1.5,0) node[midway, above] {$\mathsmaller{1}$}; 

\draw[postaction={decorate}] (1.5,1) --++ (0,-1) node[midway, right] {$\mathsmaller{2}$}; 

\draw[postaction={decorate}] (0,-1) -- (0,0) node[midway,left] {$\mathsmaller{1}$}; 

\draw[postaction={decorate}] (0,-1) -- (1.5,-1) node[midway,below] {$\mathsmaller{1}$};

\draw[postaction={decorate}] (1.5,0) --++ (0,-1) node[midway,right] {$\mathsmaller{1}$};

\draw[postaction={decorate}] (0,-2) -- (0,-1) node[midway,left] {$\mathsmaller{2}$}; 

\draw[postaction={decorate}] (1.5,-1) -- (1.5,-2) node[midway,right] {$\mathsmaller{2}$};
\end{scope}

\end{tikzpicture}
    \end{aligned}
\end{equation}
The left square can be simplified using the relation Figure \ref{MOYe}. The other two are obtained from this one by a cycle move, so by Lemma \ref{lem:cycle move}  all square faces can be simplified. Applying these reductions and removing closed innermost loops as necessary, we reduce $\brak{\Gamma}$ to a direct sum of empty webs. 

The theorem is proven by a double induction, following the proof of \cite[Theorem 2.4]{Wu} in the decategorified setting. Suppose that the statement of the theorem holds for annular webs all of whose edges and closed loops are of thickness at most $m$, for some $m\geq 2$. Additionally assume that that, for some $k\geq 0$, that the statement of the theorem holds for all webs where exactly $k$ edges and loops have thickness $m+1$. We will prove that the theorem then holds for all annular webs such that (1) all edges and loops have thickness at most $m+1$ and (2) there are exactly $k+1$ edges and loops of thickness $m+1$. Let $\Gamma$ be such a web.

Suppose a closed loop of $\Gamma$ has thickness $m+1$. It is either innermost, in which case it can be removed, or it contains another annular web, which can be simplified by using the inductive step. Suppose then that $\Gamma$ has an edge $e$ of thickness $m+1$, which is necessarily of the form shown in the left-most web in \eqref{eq:inductive webs}

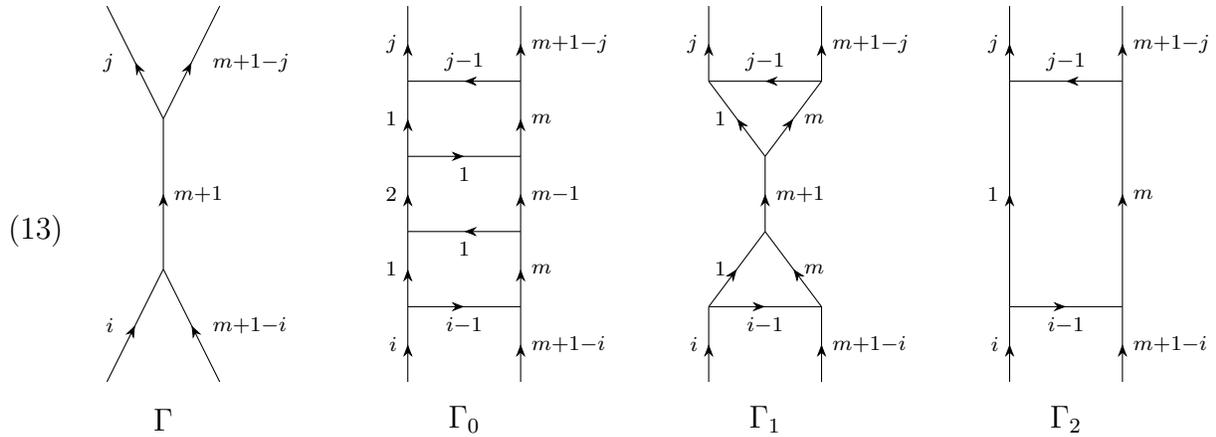
\begin{equation}
    \label{eq:inductive webs}
    \begin{aligned}
    \begin{tikzpicture}[yscale = 1,decoration={
    markings,
    mark=at position 0.5 with {\arrow{Stealth}}}
    ] 

\draw[postaction={decorate}] (0,-4) --++ (.75,1.5) node[midway,left]  {$\mathsmaller{i}\ $} ;

\draw[postaction={decorate}] (1.5,-4) --++ (-.75,1.5) node[midway,right] {$\ \mathsmaller{m+1-i}$};

\draw[postaction={decorate}] (.75,-2.5) --++ (0,2) node[midway,right] {$\mathsmaller{m+1}$};

\draw[postaction={decorate}] (.75,-.5) --++ (-.75,1.5) node[midway,left] {$\mathsmaller{j}\ $};

\draw[postaction={decorate}] (.75,-.5) --++ (.75,1.5) node[midway,right] {$\ \mathsmaller{m+1-j}$};

\node at (.75,-4.5) {$\Gamma$};

\begin{scope}[shift={(4,0)}]

\draw[postaction={decorate}] (0,0) --++ (0,1) node[midway, left] {$\mathsmaller{j}$};

\draw[postaction={decorate}] (1.5,0) --++ (0,1) node[midway, right] {$\mathsmaller{m+1-j}$}; 


\draw[postaction={decorate}] (0,-1) --++ (0,1) node[midway,left] {$\mathsmaller{1}$}; 

\draw[postaction={decorate}] (1.5,-1) --++ (0,1) node[midway,right] {$\mathsmaller{m}$};


\draw[postaction={decorate}] (0,-2) --++ (0,1) node[midway,left] {$\mathsmaller{2}$}; 

\draw[postaction={decorate}] (1.5,-2) --++ (0,1) node[midway,right] {$\mathsmaller{m-1}$};


\draw[postaction={decorate}] (0,-3) --++ (0,1) node[midway,left] {$\mathsmaller{1}$}; 

\draw[postaction={decorate}] (1.5,-3) --++ (0,1) node[midway,right] {$\mathsmaller{m}$};

\draw[postaction={decorate}] (0,-4) --++ (0,1) node[midway,left] {$\mathsmaller{i}$}; 

\draw[postaction={decorate}] (1.5,-4) --++ (0,1) node[midway,right] {$\mathsmaller{m+1-i}$};


\draw[postaction={decorate}] (1.5,0) --++ (-1.5,0) node[midway, above] {$\mathsmaller{j-1}$};

\draw[postaction={decorate}] (0,-1) --++ (1.5,0) node[midway,below] {$\mathsmaller{1}$};

\draw[postaction={decorate}] (1.5,-2) --++ (-1.5,0) node[midway,below] {$\mathsmaller{1}$};

\draw[postaction={decorate}] (0,-3) --++ (1.5,0) node[midway,below] {$\mathsmaller{i-1}$};

\node at (.75,-4.5) {$\Gamma_0$}; 
\end{scope}


\begin{scope}[shift={(8,0)}]
\draw[postaction={decorate}] (0,0) --++ (0,1) node[midway, left] {$\mathsmaller{j}$};

\draw[postaction={decorate}] (1.5,0) --++ (0,1) node[midway, right] {$\mathsmaller{m+1-j}$}; 


\draw[postaction={decorate}] (.75,-1) --++ (-.75,1) node[midway,left] {$\mathsmaller{1}$}; 

\draw[postaction={decorate}] (.75,-1) --++ (.75,1) node[midway,right] {$\mathsmaller{m}$};


\draw[postaction={decorate}] (.75,-2) --++ (0,1) node[midway,right] {$\mathsmaller{m+1}$};


\draw[postaction={decorate}] (0,-3) --++ (.75,1) node[midway,left] {$\mathsmaller{1}$}; 

\draw[postaction={decorate}] (1.5,-3) --++ (-.75,1) node[midway,right] {$\mathsmaller{m}$};

\draw[postaction={decorate}] (0,-4) --++ (0,1) node[midway,left] {$\mathsmaller{i}$}; 

\draw[postaction={decorate}] (1.5,-4) --++ (0,1) node[midway,right] {$\mathsmaller{m+1-i}$};


\draw[postaction={decorate}] (1.5,0) --++ (-1.5,0) node[midway, above] {$\mathsmaller{j-1}$};

\draw[postaction={decorate}] (0,-3) --++ (1.5,0) node[midway,below] {$\mathsmaller{i-1}$};

\node at (.75,-4.5) {$\Gamma_1$};
\end{scope}


\begin{scope}[shift={(12,0)}]
\draw[postaction={decorate}] (0,0) --++ (0,1) node[midway, left] {$\mathsmaller{j}$};

\draw[postaction={decorate}] (1.5,0) --++ (0,1) node[midway, right] {$\mathsmaller{m+1-j}$}; 


\draw[postaction={decorate}] (0,-3) --++ (0,3) node[midway,left] {$\mathsmaller{1}$}; 

\draw[postaction={decorate}] (1.5,-3) --++ (0,3) node[midway,right] {$\mathsmaller{m}$};

\draw[postaction={decorate}] (0,-4) --++ (0,1) node[midway,left] {$\mathsmaller{i}$}; 

\draw[postaction={decorate}] (1.5,-4) --++ (0,1) node[midway,right] {$\mathsmaller{m+1-i}$};


\draw[postaction={decorate}] (1.5,0) --++ (-1.5,0) node[midway, above] {$\mathsmaller{j-1}$};

\draw[postaction={decorate}] (0,-3) --++ (1.5,0) node[midway,below] {$\mathsmaller{i-1}$};

\node at (.75,-4.5) {$\Gamma_2$}; 
\end{scope}

\end{tikzpicture}
    \end{aligned}
\end{equation}

Consider the three  webs $\Gamma_0, \Gamma_1, \Gamma_2$, shown in \eqref{eq:inductive webs},  obtained by locally modifying $\Gamma$ near $e$. The categorified MOY relations imply that 
\[
\brak{\Gamma} \{ ([j]_q[i]_q,0)\}  \cong \brak{\Gamma_1},\ \ \
\brak{\Gamma_0}  \cong \brak{\Gamma_1} \oplus \brak{\Gamma_2} \{([m-1]_q,0)\}.
\]
By the inductive step, the theorem holds for $\Gamma_0$ and $\Gamma_2$. Then $\brak{\Gamma_1}$, being a direct summand of a finitely generated graded free $R_N$-module is itself free. Similarly, it then follows that $\brak{\Gamma}$ is a finitely generated and free $\RN$-module. That the graded rank can be computed using the MOY relations together with \eqref{eq:annular MOY relation} is clear from the the above arguments.
\end{proof}



\subsection{Equivariant annular \texorpdfstring{$\mathfrak{gl}_N$}{gln} link homology}
\label{sec:equivariant annular gln link homology} 
In this section we define the   chain complex associated to an exterior colored, oriented annular link $L$, following \cite[Definition 3.3]{ETW}. Let $\A = S^1 \times [0,1]$ denote the annulus. We identify the interior of $\A$ with the punctured plane $\P$. A link in the thickened annulus $\A \times [0,1]$ is then viewed as a link in $\P\times [0,1]$, and we fix a diagram $D \subset \P$ of $L$. A \emph{colored} annular link is an oriented link $L\subset \P \times [0,1]$ with each component labeled by an integer in $\{0,\ldots, N\}$; components of $D$ are then similarly labeled.

\begin{definition}
To a positive crossing with overstrand colored $i$ and understrand colored $j$, with $i\geq j$ (Figure \ref{fig:gln complex crossing}) assign the complex 
\[
\brak{\Gamma_0} \{(c_0,0)\} \to \brak{\Gamma_1} \{(c_1,0)\} \to \cdots \to \brak{\Gamma_j}\{(c_j,0)\},
\]
where the webs $\Gamma_k$  are shown in Figure \ref{fig:gln complex Gamma_k}, $\brak{\Gamma_k}$ is in homological degree $k$, the grading shifts are given by $c_k = -k - j(N-j)$ (applied only to the quantum grading), and the maps are induced by a single foam, shown in Figure \ref{fig:gln complex differential}. 

If $i<j$ or the crossing is negative, then the complex assigned to the crossing is obtained from the one above in the manner described in \cite[Definition 3.3]{ETW}. 

Finally, the chain complex $C_N(D)$ of an annular link diagram $D$ is obtained by replacing each crossing with the corresponding complex and tensoring them together in a planar algebra fashion. We let $H_N(D)$ denote the homology of this chain complex.
\end{definition}

\begin{figure}
\centering
\subcaptionbox{A positive crossing, with $i\geq j$.
\label{fig:gln complex crossing}}[.4\linewidth]
{\begin{tikzpicture}[yscale=.9]

\begin{scope}[decoration={
    markings,
    mark=at position 0.75 with {\arrow{Stealth}}}
    ]

\draw[postaction={decorate}] (1,-2) -- (-1,1);

\draw[postaction={decorate}, preaction={draw, line width=8pt, white}] (-1,-2) -- (1,1) ; 

\node[right] at (1,1) {$i$}; 
\node[left] at (-1,1) {$j$}; 
\end{scope}

\end{tikzpicture}
}
\subcaptionbox{The web $\Gamma_k$, $0 \leq k \leq j$, appearing in homological degree $k$.
\label{fig:gln complex Gamma_k}}[.4\linewidth]
{\begin{tikzpicture}[scale=1]

\begin{scope}[shift={(4,0)},  decoration={
    markings,
    mark=at position 0.75 with {\arrow{Stealth}}}
    ]

\draw[postaction={decorate}] (0,0) --++ (0,1) node[midway, left] {$\mathsmaller{j}$}; 

\draw[postaction={decorate}] (0,0) --++ (1.5,0) node[midway, above] {$\mathsmaller{i-j+k}$}; 

\draw[postaction={decorate}] (1.5,0) --++ (0,1) node[midway, right] {$\mathsmaller{i}$}; 

\draw[postaction={decorate}] (0,-1) --++ (0,1) node[midway,left] {$\mathsmaller{i+k}$}; 

\draw[postaction={decorate}] (1.5,-1) --++ (-1.5,0) node[midway,below] {$\mathsmaller{k}$};

\draw[postaction={decorate}] (1.5,-1) --++ (0,1) node[midway,right] {$\mathsmaller{j-k}$};

\draw[postaction={decorate}] (0,-2) --++ (0,1) node[midway,left] {$\mathsmaller{i}$}; 

\draw[postaction={decorate}] (1.5,-2) --++ (0,1) node[midway,right] {$\mathsmaller{j}$};

\end{scope}

\end{tikzpicture}
}\\ \vskip2ex
\subcaptionbox{The foam cobordism $\Gamma_k \to \Gamma_{k+1}$ inducing the differential. The shaded facet has thickness $1$. Thicknesses and orientations of facets are determined by $\Gamma_k$ and $\Gamma_{k+1}$. 
\label{fig:gln complex differential}}[1\linewidth]
{\begin{tikzpicture}[scale=.8]

\begin{scope}[shift={(15,0)},  decoration={
    markings,
    mark=at position 0.4 with {\arrow{Stealth}}}
    ]
    \draw[postaction={decorate}] (0,-2) -- (-2,-2);
    \draw[postaction={decorate}] (2,-2) -- (0,-2);
      \draw[postaction={decorate}] (4,-2) -- (2,-2);
      
        \draw[postaction={decorate}] (0,-2) -- (0,0);
        \draw[postaction={decorate}] (2,0) -- (2,-2);
        
        \draw (-2,-2) -- (-2,2);
        \draw (4,-2) -- (4,2);
        \draw (5,-1) -- (5,3);
        \draw[densely dashed] (-1,-1) -- (-1,2);
        \draw (-1,2) -- (-1,3); 
        
        \draw[postaction={decorate}] (2,0) -- (0,0);
        \draw[postaction={decorate}] (0,0) -- (0,2);
        \draw[postaction={decorate}] (2,2) -- (2,0);
        
        \draw[densely dashed,postaction={decorate}] (0,0) -- (1,1); 
        \draw[densely dashed,postaction={decorate}] (3,1) -- (2,0); 
        \draw[densely dashed,postaction={decorate}] (1,1) -- (3,1); 
        \fill[gray,opacity=0.40] (0,0) -- (1,1) -- (3,1) -- (2,0);
        
        \draw[densely dashed,postaction={decorate}] (1,2) -- (1,1) ; 
        \draw[densely dashed] (3,1) -- (3,2); 
        
        \draw (1,2) -- (1,3);
        \draw[postaction={decorate}] (3,2) -- (3,3); 
        \draw[postaction=decorate] (0,2) -- (1,3); 
        
        \draw[densely dashed, postaction={decorate}] (1,1) -- (1,-1); 
        \draw[densely dashed, postaction={decorate}] (3,-1) -- (3,1);
        \draw[densely dashed, postaction={decorate}] (0,-2) -- (1,-1);
        \draw[densely dashed, postaction={decorate}] (3,-1) -- (2,-2);
        \draw[densely dashed, postaction={decorate}] (3,-1) -- (1,-1);
        \draw[densely dashed, postaction={decorate}] (1,-1) -- (-1,-1);
        \draw[densely dashed] (4,-1) -- (3,-1);
        \draw[ postaction={decorate}] (5,-1) -- (4,-1);
        
        \draw[ postaction={decorate}] (0,2) -- (-2,2);
        \draw[ postaction={decorate}] (1,3) -- (-1,3); 
        \draw[ postaction={decorate}] (3,3) -- (1,3); 
        \draw[ postaction={decorate}] (5,3) -- (3,3);
        \draw[ postaction={decorate}] (2,2) -- (0,2); 
        \draw[ postaction={decorate}] (4,2) -- (2,2);
        
        \draw[ postaction={decorate}] (3,3) -- (2,2);
        
        \node at (6,-2) {$\Gamma_k$};
        \node at (6,3) {$\Gamma_{k+1}$};
        
        \node at (1.5,.5) {$\mathsmaller{1}$};
\end{scope}

\end{tikzpicture}
}
\caption{}\label{fig:gln complex}
\end{figure}

The differential is induced by foams which are disjoint from $\L$ and hence preserves annular degree. Chain groups of $C_N(D)$ and homology groups $H_N(D)$ are $\Z\oplus \Z^N$-graded $\RN$-modules. 

Consider colored and oriented annular link diagrams $D_0$ and $D_1$ where $D_1$ is obtained from $D_0$ by  a Reidemeister move or a Morse move (cup, cap, or saddle) which occurs away from the puncture; saddles must involve strands of the same color. There is an induced chain map $C_N(D_0) \to C_N(D_1)$, given by the natural foam map for Morse moves and otherwise given by chain homotopy equivalences for Reidemeister moves. This map has bidegree $(d,0, \ldots, 0)$, where $d=0$ for Reidemeister moves,  $d = -i(N-i)$ for a cup or a cap involving a circle of color $i$, and $d= i(N-i)$ for a saddle involving strands of color $i$. 

A \emph{colored} cobordism is an oriented link cobordism $S\subset \RR^3\times [0,1]$ in which each component is labeled by an element of $\{0, \ldots, N\}$. The boundary links $S\cap \RR^3\times \{0\}$ and $S\cap \RR^3 \times \{1\}$ are then naturally colored. Any colored link cobordism can be represented as a sequence of the elementary Reidemeister or Morse cobordisms described above. A \emph{cobordism} between colored oriented links $L, L'$ is a colored cobordism $S$ such that the induced coloring on the boundary of $S$ agrees with the coloring of $L$ and  $L'$. These notions extend in a straightforward manner to cobordisms in $\P\times [0,1] \times [0,1]$ between annular links. 

\begin{proposition}
Let $L$ be an annular link with a diagram $D$. If $D'$ is an annular link diagram obtained from $D$ by a Reidemeister move away from the puncture, then $C_N(D)$ and $C_N(D')$ are chain homotopy equivalent. Consequently, $H_N(D)$ is an invariant of $L$ up to isomorphism.

Let $S \subset \A\times [0,1] \times [0,1]$ be a colored cobordism between annular links $L$ and $L'$. Pick diagrams $D, D'$ for $L, L'$ and represent $S$ as a sequence of elementary cobordisms. Up to chain homotopy equivalence, the resulting chain map $C_N(D) \to C_N(D')$ is independent of the choice of decomposition of $S$ into elementary pieces. 
\end{proposition}

\begin{proof}
Local invariance of colored $\gl_N$ homology is established in  \cite[Theorem 3.5]{ETW}. Since annular link diagrams representing isotopic annular links are related by Reidemeister moves away from the puncture, it follows that $C_N(D)$ and $C_N(D')$ are chain homotopy equivalent. 

Isotopic annular cobordisms can be related by movie moves occuring away from the puncture, so the statement on functoriality follows from \cite[Theorem 4.5]{ETW}.
\end{proof}

\bibliographystyle{amsalpha}
\bibliography{main}

\end{document}